\documentclass[12pt]{article}
\newcommand{\real}{{\mathbb R}}
\newcommand{\complex}{{\mathbb C}}
\newcommand{\intplus}{{\mathbb N}}
\newcommand{\allint}{{\mathbb Z}}
\renewcommand{\div}{\mathop{\rm div}}
\newcommand{\rot}{\mathop{\rm rot}}

\renewcommand{\span}{{\rm span}}
\newcommand{\range}{{\rm range}}
\newcommand{\e}{{\rm e}}		%% \e = 2.771828...
\renewcommand{\d}{\,{\rm d}}		%% integration measure
\def\I{{\rm i}}				%% \I = \sqrt{-1}

\newcommand{\const}{{\rm const}}
\newcommand{\sign}{\mathop{\rm sign}}

\newcommand{\cE}{{\cal E}}
\newcommand{\cL}{{\cal L}}
\newcommand{\cO}{{\cal O}}
\newcommand{\cV}{{\cal V}}
\newcommand{\cW}{{\cal W}}

\newcommand{\cc}{{\bf c}}
\newcommand{\ee}{{\bf e}}
\newcommand{\ff}{{\bf f}}
\renewcommand{\gg}{{\bf g}}
\newcommand{\hh}{{\bf h}}
\newcommand{\pp}{{\bf p}}
\newcommand{\qq}{{\bf q}}
\newcommand{\rr}{{\bf r}}
\newcommand{\uu}{{\bf u}}
\newcommand{\vv}{{\bf v}}
\newcommand{\ww}{{\bf w}}
\renewcommand{\AA}{{\bf A}}
\newcommand{\BB}{{\bf B}}
\newcommand{\FF}{{\bf F}}
\newcommand{\RR}{{\bf R}}
\newcommand{\UU}{{\bf U}}
\newcommand{\oomega}{{\mbox{\boldmath$\omega$}}}
\newcommand{\PPhi}{{\mbox{\boldmath$\Phi$}}}
\newcommand{\oone}{{\bf 1}}

\renewcommand{\L}{{\rm L}}

\newtheorem{theorem}{Theorem}[section]
\newtheorem{lemma}[theorem]{Lemma}

\newtheorem{proposition}[theorem]{Proposition}
\newtheorem{corollary}[theorem]{Corollary}
\newtheorem{remark}[theorem]{Remark}

\newcommand{\reff}[1]{(\ref{#1})}

\newcommand{\proof}{{\noindent \bf Proof:\ }}
\newcommand{\half}{{\frac{1}{2}}}

\def\build#1_#2^#3{\mathrel{
  \mathop{\kern 0pt#1}\limits_{#2}^{#3}}}
\def\QED{\mbox{}\hfill QED}

\newdimen\texpscorrection
\newdimen\figcenter
\def\figurewithtex #1 #2 #3 #4 #5\cr{\null
  {\goodbreak\figcenter=\hsize\relax
  \advance\figcenter by -#4truecm
  \divide\figcenter by 2
  \begin{figure}[hbt]
  \vskip #3truecm\noindent\hskip\figcenter
  \includegraphics{#1}{\hskip\texpscorrection\input #2 }
  \vskip 0.8truecm\noindent \vbox{\noindent {\footnotesize #5}}
  \end{figure}}}
\def\point#1 #2 #3 {\rlap{\kern #1 truecm
\raise #2 truecm \hbox{#3}}}

\usepackage{amsfonts}

\begin{document}

\title{Long-time asymptotics of the Navier-Stokes and 
vorticity equations on $\real^3$}

\author{Thierry Gallay \\ 
Universit\'e de Paris-Sud \\ 
Math\'ematiques \\
B\^atiment 425 \\
F-91405 Orsay \\
France
\and
C. Eugene Wayne \\
Department of Mathematics \\
and  Center for BioDynamics \\
Boston University \\
111 Cummington Street\\
Boston, MA 02215, USA}
\date{February 23, 2001}
\maketitle

\begin{abstract}
We use the vorticity formulation to study the long-time behavior of
solutions to the Navier-Stokes equation on $\real^3$. We assume that
the initial vorticity is small and decays algebraically at infinity. 
After introducing self-similar variables, we compute the long-time 
asymptotics of the rescaled vorticity equation up to second order. 
Each term in the asymptotics is a self-similar divergence-free vector 
field with Gaussian decay at infinity, and the coefficients in 
the expansion can be determined by solving a finite system of 
ordinary differential equations. As a consequence of our results, 
we are able to characterize the set of solutions for which the 
velocity field satisfies $\|\uu(\cdot,t)\|_{L^2} = o(t^{-5/4})$ as
$t \to +\infty$. In particular, we show that these solutions lie 
on a smooth invariant submanifold of codimension $11$ in our 
function space.
\end{abstract}

%%%%%%%%%%%%%%%%%%%%%%%%%%%%%%%%%%%%%%%%%%%%%%%%%%%%%%%%%%%%%%%%%%%%%%%%%%%

\section{Introduction}\label{intro}

We consider the motion of an incompressible viscous fluid filling
the whole space $\real^3$. If no external force is applied, the 
velocity $\uu(x,t)$ of the fluid satisfies the Navier-Stokes equation
\begin{equation}\label{NS}
  \partial_t \uu + (\uu\cdot\nabla)\uu = \nu \Delta \uu - \frac{1}{\rho}
  \nabla p\ , \quad \div \uu = 0\ ,
\end{equation}
where $\rho$ is the density of the fluid, $\nu$ is the kinematic viscosity, 
and $p(x,t)$ is the pressure field. Replacing $x,t,\uu,p$ with the 
dimensionless quantities
$$
  \frac{x}{L}\ , \quad \frac{\nu t}{L^2}\ , \quad \frac{L \uu}{\nu}\ ,
  \quad \frac{L^2 p}{\rho \nu^2}\ ,
$$
where $L$ is an arbitrary length scale, Eq.\reff{NS} is transformed 
into
\begin{equation}\label{NS3}
  \partial_t \uu + (\uu\cdot\nabla)\uu = \Delta \uu - \nabla p\ , \quad
  \div \uu = 0\ .
\end{equation}
Since the length $L$ was arbitrary, Eq.\reff{NS3} is still invariant 
under the scaling transformation
\begin{equation}\label{scaltransf}
  \uu(x,t) \mapsto \lambda \uu(\lambda x,\lambda^2 t)\ , \quad
  p(x,t) \mapsto \lambda^2 p(\lambda x,\lambda^2 t)\ ,
\end{equation}
for any $\lambda > 0$. 

As no external force is applied, it is intuitively clear that all 
finite-energy solutions of \reff{NS3} should converge, as time goes
to infinity, to the rest state $\uu \equiv 0$, $p \equiv \const$. As a
matter of fact, if $\uu(x,t)$ is any global weak solution in $L^2(\real^3)$
satisfying the energy inequality, it is known that $\|\uu(\cdot,t)\|_{L^2}
\to 0$ as $t \to \infty$ (see Wiegner \cite{wiegner:1987}). Moreover, if 
\begin{equation}\label{lincond}
  \|\e^{t\Delta}\uu(\cdot,0)\|_{L^2} \le \frac{C}{(1+t)^\alpha}\ , \quad
  t \ge 0\ ,
\end{equation}
for some $\alpha \ge 0$, then
\begin{equation}\label{nlincond}
  \|\uu(\cdot,t)\|_{L^2} \le \frac{C'}{(1+t)^\beta}\ , \quad
  t \ge 0\ ,
\end{equation}
where $\beta = \min(\alpha,5/4)$. This result shows that the solutions
of \reff{NS3} decay to zero at the same rate as those of the linear heat 
equation, provided this rate does not exceed $t^{-5/4}$. As we 
shall see below, the restriction $\beta \le 5/4$ in \reff{nlincond} is
due to the nonlinearity in \reff{NS3} and to the incompressibility 
condition $\div \uu = 0$. 

Wiegner's result raises a very natural question: can we characterize the 
set of solutions of \reff{NS3} such that $t^{5/4}\|\uu(\cdot,t)\|_{L^2}
\to 0$ as $t \to \infty$? Put differently, given a solution $\uu(x,t)$ 
satisfying \reff{nlincond} with $\beta = 5/4$, under which conditions
can we prove the corresponding {\it lower bound} $\|\uu(\cdot,t)\|_{L^2} 
\ge C (1+t)^{-5/4}$? 
This problem has been intensively studied during the last 15 years, 
especially by M.E.~Schonbek \cite{schonbek:1985}, \cite{schonbek:1986}, 
\cite{schonbek:1991}, \cite{schonbek:1992}, who found sufficient 
conditions for such a lower bound to hold. For technical reasons, 
these results were established assuming some additional decay of the 
initial data $\uu_0 = \uu(\cdot,0)$. Typically, it was assumed that 
$\uu_0 \in L^2(\real^3)^3$ and $(1+|x|)\uu_0 \in L^1(\real^3)^3$, 
so that \reff{lincond} holds with $\alpha = 5/4$. 

Very recently, T.~Miyakawa and M.E.~Schonbek obtained an interesting 
characterization of the ``rapidly decreasing'' solutions of the 
Navier-Stokes equation in $\real^N$, $N \ge 2$. In the case $N = 3$, 
their result reads:

\begin{theorem}\label{MiSchon} {\bf \cite{miyakawa:2000}}
Assume that $\uu_0 \in L^2(\real^3)^3$, $\div \uu_0 = 0$, and
$(1+|x|)\uu_0 \in L^1(\real^3)^3$. Let $\uu(x,t)$ be a global
weak solution of \reff{NS3} with initial data $\uu(\cdot,0) = \uu_0$,
satisfying the bound \reff{nlincond} with $\beta = 5/4$. For all 
$k,\ell \in \{1,2,3\}$, define
\begin{equation}\label{bcmom}
  b_{k \ell} = \int_{\real^3} x_k u_\ell(x,0) \d x\ , \quad
  c_{k \ell} = \int_0^{\infty} \int_{\real^3} u_k(x,t) u_{\ell}(x,t) 
  \d x \d t\ .
\end{equation}
Then 
\begin{equation}\label{optidecay}
  \lim_{t\to\infty} t^{5/4}\|\uu(\cdot,t)\|_{L^2} = 0
\end{equation}
if and only if there exists $c \ge 0$ such that
\begin{equation}\label{more-moments}
  b_{k \ell} = 0 \quad \mbox{and }\quad c_{k \ell} = c \delta_{k \ell}\ ,
  \quad k,\ell \in \{1,2,3\}\ .
\end{equation}  
\end{theorem}

The proof is a direct calculation using the integral equation satisfied
by the solutions of \reff{NS3}. As the authors themselves remark, 
this argument does not provide much intuition as to the meaning of the 
conditions \reff{more-moments}. From our point of view, the most 
surprising feature of Theorem~\ref{MiSchon} is the fact that assertion
\reff{optidecay} is {\it translation invariant} in time, whereas
conditions \reff{more-moments} are not. More precisely, if a solution
$\uu(x,t)$ satisfies \reff{optidecay}, so does $\uu(x,t+T)$ for any
$T > 0$. However, the property that $(1+|x|)\uu$ be integrable is not
preserved under the Navier-Stokes evolution (even for strong solutions, 
where uniqueness holds). Thus, if $(1+|x|)\uu_0 \in L^1(\real^3)^3$, 
then in general $(1+|x|)\uu(\cdot,T) \notin L^1(\real^3)^3$ for $T > 0$, 
so that the moments $b_{k\ell}$ do not make sense if $\uu(x,0)$ is 
replaced by $\uu(x,T)$. Moreover, if the matrix $(c_{k\ell})$ is scalar, 
this property will be lost if $\uu(x,t)$ is replaced by $\uu(x,t+T)$. 
Summarizing, Theorem~\ref{MiSchon} is a characterization of those 
rapidly decreasing solutions of \reff{NS3} {\it whose initial data
lie in the noninvariant subspace} $\{\uu_0 \in L^2(\real^3)^3\,|\, 
(1+|x|)\uu_0 \in L^1(\real^3)^3\}$. Nontrivial examples of such 
solutions have been recently constructed by L.~Brandolese 
\cite{brandolese:2001}. 

In this paper, we use the vorticity formulation to study the 
long-time behavior of the solutions of the Navier-Stokes equation
\reff{NS3}. Setting $\oomega = \rot \uu$, Eq.\reff{NS3} is transformed
into
\begin{equation}\label{V3}
  \partial_t \oomega +(\uu \cdot \nabla) \oomega -(\oomega \cdot \nabla ) 
  \uu = \Delta \oomega\ , \quad \div \oomega = 0\ .
\end{equation}
The velocity field $\uu$ can be reconstructed from $\oomega$ via the 
Biot-Savart law:
\begin{equation} \label{BS3}
  \uu(x) = -\frac{1}{4\pi} \int_{\real^3} 
  \frac{(x-y) \wedge \oomega(y)}{|x - y|^3} \d y\ ,
  \quad x \in \real^3\ ,
\end{equation}
where $\wedge$ denotes the cross product in $\real^3$. Although \reff{NS3}
and \reff{V3} are formally equivalent, we believe that using the vorticity
formulation to compute the long-time asymptotics has a crucial advantage:
roughly speaking, the spatial decay of $\oomega$ is preserved under the 
evolution defined by \reff{V3}. For instance, if $(1+|x|)^m \oomega_0 \in
L^2(\real^3)^3$ for some $m \ge 0$, then \reff{V3} has a unique local 
solution $\oomega(x,t)$ with initial data $\oomega_0$ satisfying 
$(1+|x|)^m \oomega(\cdot,t) \in L^2(\real^3)^3$ whenever it exists. 
Again, we point out that this property does {\it not} hold for the 
velocity field $\uu(x,t)$ if $m \ge 5/2$. This is the reason why the 
integrability condition $(1+|x|)\uu \in L^1(\real^3)^3$ is not preserved
under evolution. 

In the sequel, we always assume that the vorticity $\oomega(x,t)$ 
is small and decreases sufficiently fast as $|x| \to \infty$. The smallness
assumption is not a restriction as far as the long-time behavior is 
concerned, since all global solutions of the Navier-Stokes equation 
(in the energy space) converge to zero as $t \to \infty$. Moreover, this
hypothesis allows to deal with global strong solutions of \reff{V3}. 
On the other hand, assuming that the vorticity decreases rapidly as 
$|x| \to \infty$ is very reasonable from a physical point of view. This
is the case, for instance, if the initial data are created by stirring the 
fluid with a (finite size) tool. In addition, this property is very
helpful to study the long-time asymptotics, since in parabolic equations
the spatial and temporal behaviors of the solutions are intimately 
connected. 

To actually compute the asymptotics, we express the vorticity 
$\oomega(x,t)$ in terms of the self-similar variables $(\xi,\tau)$ 
defined by $\xi = x/\sqrt{1{+}t}$, $\tau = \log(1{+}t)$, see 
\reff{omega-w3} below. Although the transformation is time-dependent, 
the rescaled vorticity $\ww(\xi,\tau)$ still satisfies an autonomous
equation, as a consequence of the scaling invariance \reff{scaltransf}. 
Linearizing this equation around the origin $\ww = 0$, we find that
the generator $\Lambda$ of the time evolution has a countable set of real, 
isolated eigenvalues with finite multiplicities, and that the essential
spectrum can be pushed arbitrarily far away into the left-half plane
by choosing the function space (i.e., the spatial decay of the vorticity)
appropriately. Thus, the long-time asymptotics in a neighborhood of the
origin are determined, at any prescribed order, by a finite system of
ordinary differential equations. 

This reduction procedure, or some variant of it, has been often applied
to investigate the long-time behavior of solutions of nonlinear parabolic
or damped hyperbolic equations \cite{bricmont:1996}, 
\cite{eckmann:1997}, \cite{eckmann:1998}, \cite{escobedo:1995},
\cite{galaktionov:1991}, \cite{gallay:1998}, \cite{gallay:2000},
\cite{kavian:1987}, \cite{wayne:1997}, \cite{zuazua:1993}. In the 
context of the
Navier-Stokes equation, rescaling techniques were used by A.~Carpio 
\cite{carpio:1994}, \cite{carpio:1996} to study the vorticity equations 
in two and three dimensions. In \cite{cannone:1996}, M.~Cannone and 
F.~Planchon constructed a large family of self-similar solutions of the
three-dimensional Navier-Stokes equation. These solutions correspond to 
fixed points of our rescaled vorticity equation, but do not belong to 
the function spaces we use, because they decay too slowly as $|x| \to 
\infty$. In a companion paper \cite{gallay:2001}, we follow the procedure
outlined above to study the solutions of the two-dimensional Navier-Stokes
and vorticity equations. In addition, we exploit the fact that the 
spectrum of the generator $\Lambda$ is discrete to construct 
finite-dimensional {\it invariant manifolds} that are approached, at 
a prescribed rate, by all solutions in a neighborhood of the origin. 

The rest of this paper is organized as follows. In Section~\ref{cauchy}, 
we prove the existence of global solutions of the vorticity equation
\reff{V3} in a neighborhood of the origin, and we estimate their decay
rate as $t \to \infty$. The results we obtain are comparable to 
those of M.~Wiegner in \cite{wiegner:1987}. Section~\ref{first} is 
devoted to the first order asymptotics. Under appropriate conditions, 
we show that
$$
   \oomega(x,t) \sim \sum_{i=1}^3 \frac{b_i}{(1+t)^2} \,\ff_i\Bigl(
   \frac{x}{\sqrt{1+t}}\Bigr)\ , \quad t \to \infty\ ,
$$
where $\ff_1, \ff_2, \ff_3$ are explicit divergence-free vector fields
with Gaussian decay at infinity, and $b_1, b_2, b_3$ are real coefficients
which can be computed from the initial data. Using \reff{BS3}, a similar
result can be obtained for the velocity field $\uu(x,t)$. In 
Section~\ref{second}, we give a higher order asymptotic expansion 
of $\oomega(x,t)$, including terms of the form $(1+t)^{-5/2}
\gg(x/\sqrt{1+t})$. This result is used in Section~\ref{ssmanif} to 
characterize the set of solutions $\oomega(x,t)$ of \reff{V3} for which
the velocity field $\uu(x,t)$ satisfies \reff{optidecay}. It is shown 
that these solutions lie on a smooth invariant manifold of finite 
codimension, which is tangent at the origin to a spectral subspace 
of the generator $\Lambda$. Intersecting this manifold with the 
(noninvariant) subspace $\{\uu_0 \in L^2(\real^3)^3\,|\, (1+|x|)\uu_0 
\in L^1(\real^3)^3\}$, we recover exactly conditions \reff{more-moments}
in Theorem~\ref{MiSchon}. Finally, Appendix~\ref{spectrum} describes 
the spectral properties of the generator $\Lambda$, and 
Appendix~\ref{velocity3d} collects various estimates of the velocity 
field $\uu$ in terms of the vorticity $\oomega$ in weighted Lebesgue
spaces. 

\medskip\noindent{\bf Notation.} Throughout the paper, 
we use boldface letters for vector-valued functions, such as 
$\uu(x,t)$ and $\oomega(x,t)$. However, to avoid a proliferation
of boldface symbols, we use standard italic characters for 
vector variables, such as $x = (x_1,x_2,x_3)$. In both cases, 
$|\cdot|$ denotes the Euclidean norm in $\real^3$: $|\uu| = 
(u_1^2+u_2^2+u_3^2)^{1/2}$, $|x| = (x_1^2+x_2^2+x_3^2)^{1/2}$.
For any $p \in [1,\infty]$, we denote by $|f|_p$ the norm of
a function $f$ in the Lebesgue space $L^p(\real^3)$. If $\ff \in
L^p(\real^3)^3$, we set $|\ff|_p = |\,|\ff|\,|_p$. Weighted norms 
play a very important role in this paper. We always denote by
$\rho : \real^3 \to \real$ the weight function defined by 
$\rho(x) = 1+|x|$. For any $m \ge 0$, we set $\|f\|_m = 
|\rho^m f|_2$, and $\|\ff\|_m = |\rho^m \ff|_2$. If $f \in 
C^0([0,T],L^p(\real^3))$, we often write $f(\cdot,t)$ or 
simply $f(t)$ to denote the map $x \mapsto f(x,t)$. Finally, 
we denote by $C$ a generic positive constant, which may differ
from place to place, even in the same chain of inequalities.

\medskip\noindent
{\bf Acknowledgements.} Part of this work was done when C.E.W.
visited the University of Paris-Sud and Th.G. the Department of 
Mathematics and Center for BioDynamics of Boston University. 
The hospitality of both institutions is gratefully acknowledged.
We also thank L.~Brandolese, I.~Gallagher, A.~Mielke, G.~Raugel, 
J.-C.~Saut, and M.~Vishik for stimulating discussions.  We are
especially indebted to A.~Mielke for bringing to our attention
the work of \cite{miyakawa:2000}, which triggered our interest 
in this problem. The research of C.E.W. is supported in part by 
the NSF under grant DMS-9803164. 

%%%%%%%%%%%%%%%%%%%%%%%%%%%%%%%%%%%%%%%%%%%%%%%%%%%%%%%%%%%%%%%%%%%%%%%%%%%

\section{The Cauchy problem for the vorticity equation}
\label{cauchy}

The aim of this section is to prove the existence of global solutions
of the vorticity equation for small initial data in weighted 
Lebesgue spaces. We first recall a few standard estimates for the 
velocity field $\uu$ in terms of the associated vorticity $\oomega = 
\rot \uu$. Further estimates in weighted spaces can be found in 
Appendix~\ref{velocity3d}. 

\begin{lemma}\label{HLS3} Let $\uu$ be the velocity field obtained from
$\oomega$ via the Biot-Savart law \reff{BS3}.\\
{\bf (a)} Assume that $1 < p < q < \infty$ and $\frac{1}{q} = \frac{1}{p} 
-\frac{1}{3}$. If $\oomega \in L^p(\real^3)^3$, then $\uu \in 
L^q(\real^3)^3$, and there exists $C>0$ such that
\begin{equation}\label{HLSbis}
  |\uu|_q \le C |\oomega|_p\ .
\end{equation}
{\bf (b)} Assume that $1 \le p < 3 < q \le \infty$, and define $\alpha \in 
(0,1)$ by the relation $\frac{1}{3} = \frac{\alpha}{p} + \frac{1-\alpha}{q}$. 
If $\oomega \in L^p(\real^3)^3 \cap L^q(\real^3)^3$, then $\uu \in 
L^\infty(\real^3)^3$, and there exists $C>0$ such that
\begin{equation}\label{interpol3}
  |\uu|_\infty \le C |\oomega|_p^\alpha |\oomega|_q^{1-\alpha}\ .
\end{equation}
{\bf (c)} Assume that $1 < p < \infty$. If $\oomega \in L^p(\real^3)^3$, then 
$\nabla \uu \in L^p(\real^3)^9$ and there exists $C>0$ such that
\begin{equation}\label{Calderon}
  |\nabla\uu|_p \le C |\oomega|_p\ .
\end{equation}
In addition, $\div \uu = 0$ and, if $\div\oomega = 0$, then $\rot \uu 
= \oomega$.
\end{lemma}

\proof Part {\bf (a)} is a direct consequence of \reff{BS3} and of 
the Hardy-Littlewood-Sobolev inequality, see for instance Stein 
\cite{stein:1970}, Theorem~V.1. To prove {\bf (b)}, assume that 
$\oomega \not\equiv 0$, and let $R = (|\omega|_p/|\omega|_q)^\beta$, 
where $\beta = \frac{\alpha}{1-3/q} = \frac{1-\alpha}{3/p-1}$. Using 
H\"older's inequality, we find
\begin{eqnarray*}
  |\uu(x)| &\le& \frac{1}{4\pi} \int_{|y|\le R} |\oomega(x-y)|
  \frac{1}{|y|^2}\d y + \frac{1}{4\pi} \int_{|y|\ge R} |\oomega(x-y)|
  \frac{1}{|y|^2}\d y \\
  &\le& C |\oomega|_q R^{1-\frac{3}{q}} + C|\oomega|_p 
   \frac{1}{R^{\frac{3}{p}-1}} \le 2C |\oomega|_p^\alpha 
   |\oomega|_q^{1-\alpha}\ .
\end{eqnarray*} 
Finally, $\nabla \uu$ is obtained from $\oomega$ via a singular integral 
kernel of Calder\'on-Zygmund type, hence \reff{Calderon} follows 
from Theorem~II.3 in \cite{stein:1970}. \QED

\medskip
In the sequel, for any $p \in [1,\infty]$, we denote by $\L^p(\real^3)$ 
the function space
\begin{equation}\label{boldL}
  \L^p(\real^3) = \bigl\{\ff \in L^p(\real^3)^3 \,|\, \div \ff = 0\bigr\}\ ,
\end{equation}
equipped with the same norm as $L^p(\real^3)^3$. As is well-known, the
$L^3$-norm of the velocity field $\uu(x,t)$ is invariant under the 
scaling transformation \reff{scaltransf}. For the vorticity $\oomega(x,t)$, 
the corresponding critical space is $\L^{3/2}(\real^3)$. The following 
result shows that the Cauchy problem for \reff{V3} is globally well-posed 
for small initial data in $\L^{3/2}(\real^3)$.

\begin{theorem}\label{LV3} There exists $\epsilon_0 > 0$ such that,
for all initial data $\oomega_0 \in \L^{3/2}(\real^3)$ with 
$|\oomega_0|_{3/2} \le \epsilon_0$, \reff{V3} has a unique solution 
$\oomega \in C^0([0,\infty),\L^{3/2}(\real^3)) \cap C^0((0,\infty),
\L^\infty(\real^3))$ satisfying $\oomega(0) = \oomega_0$. Moreover, for all 
$p \in [\frac{3}{2},+\infty]$, there exists $C_p > 0$ such 
that
\begin{equation}\label{lpvort}
  |\oomega(t)|_p \le \frac{C_p |\oomega_0|_{3/2}}{t^{1-\frac{3}{2p}}}\ ,
  \quad t>0\ .
\end{equation}
Finally, if $\uu(x,t)$ is the velocity field obtained from
$\oomega(x,t)$ via the Biot-Savart law \reff{BS3}, then $\uu(\cdot,t) 
\in \L^q(\real^3)$ for all $q \in [3,+\infty]$ if $t > 0$, and there exists 
$C_q > 0$ such that
\begin{equation}\label{lqvel}
  |\uu(t)|_q \le \frac{C_q |\oomega_0|_{3/2}}{t^{\half - \frac{3}{2q}}}\ ,
  \quad t>0\ .
\end{equation}
\end{theorem}

\proof The proof of Theorem~\ref{LV3} follows exactly the argument
of Kato \cite{kato:1984} which shows that the Navier-Stokes equation
has global solutions for small initial data in $\L^3(\real^3)$. 
The same argument also shows that the Cauchy problem for \reff{V3} is 
{\rm locally} well-posed in $\L^{3/2}(\real^3)$, without smallness 
assumption on the data. More generally, one can prove that \reff{V3}
has global solutions for small data in the Morrey space $M^{3/2}(\real^3)$, 
see \cite{giga:1989}.

\medskip Following \cite{gallay:2001}, we now introduce the 
``scaling variables'' 
\begin{equation}\label{scalvar}
   \xi = \frac{x}{\sqrt{1+t}}\ , \quad \tau = \log(1+t)\ .
\end{equation}
If $\oomega(x,t)$ is a solution of \reff{V3} and if $\uu(x,t)$ is
the corresponding velocity field, we set
\begin{eqnarray}\label{omega-w3}
 \oomega(x,t) &=& \frac{1}{1+t} \ \ww\Bigl(\frac{x}{\sqrt{1+t}},
 \log(1+t)\Bigr)\ , \\ \label{u-v3}
 \uu(x,t) &=& \frac{1}{\sqrt{1+t}} \ \vv\Bigl(\frac{x}{\sqrt{1+t}},
 \log(1+t)\Bigr)\ .
\end{eqnarray}
Then the rescaled vorticity $\ww(\xi,\tau)$ satisfies the 
evolution equation
\begin{equation}\label{SV3}
  \partial_\tau \ww = \Lambda \ww - ( \vv \cdot \nabla ) \ww +
  (\ww \cdot \nabla) \vv\ , \quad \div \ww = 0\ ,
\end{equation}
where $\Lambda$ is the differential operator
\begin{equation}\label{Lamop}
  \Lambda = \Delta_\xi + \frac{1}{2}\xi \cdot \nabla_\xi + 1\ ,
  \quad \xi \in \real^3\ .
\end{equation} 
The rescaled velocity $\vv$ is reconstructed from $\ww$ via the 
Biot-Savart law:
\begin{equation} \label{SBS3}
  \vv(\xi) = -\frac{1}{4\pi} \int_{\real^3} 
  \frac{(\xi-\eta) \wedge \ww(\eta)}{|\xi - \eta|^3} \d\eta\ .
\end{equation}

As in the two-dimensional case \cite{gallay:2001}, we shall solve the 
rescaled vorticity equation in weighted $L^2$ spaces. For any $m \ge 0$, 
we define
\begin{equation}\label{L2m}
  L^2(m) = \bigl\{f \in L^2(\real^3) \,|\, \|f\|_m < \infty \bigr\}\ ,
\end{equation}
where
$$
  \|f\|_m = \left(\int_{\real^3}(1+|\xi|)^{2m} |f(\xi)|^2 \d\xi 
  \right)^{1/2} = |\rho^m f|_2\ .
$$
Here and in the sequel, we denote by $\rho$ the weight function 
$\rho(\xi) = 1+|\xi|$. In analogy with \reff{boldL}, we introduce the 
space of divergence free vector fields
\begin{equation}\label{LL2m}
  \L^2(m) = \{\ff \in L^2(m)^3 \,|\, \div \ff = 0\}\ ,
\end{equation}
equipped with the norm $\|\ff\|_m = |\rho^m|\ff||_2$, where $|\ff| = 
(f_1^2+f_2^2+f_3^2)^{1/2}$.

In Appendix~\ref{spectrum}, we show that the operator $\Lambda$ 
is the generator of a strongly continuous semigroup $\e^{\tau\Lambda}$
in $\L^2(m)$, for any $m \ge 0$. Since $\partial_i\Lambda = 
(\Lambda+\frac{1}{2})\partial_i$ for $i = 1,2,3$ (where $\partial_i = 
\partial_{\xi_i}$), it is clear that $\partial_i \e^{\tau\Lambda} = 
\e^{\frac{\tau}{2}}\e^{\tau\Lambda}\partial_i$ for all $\tau \ge 0$. 
Thus, using the fact that $\div \vv = \div \ww = 0$, we can rewrite 
\reff{SV3} in integral form as follows:
\begin{equation}\label{SV3int}
  w_i(\tau) = \e^{\tau\Lambda}w_i(0) + \sum_{j=1}^3 \int_0^\tau
  \partial_j \e^{(\tau-s)(\Lambda-\frac{1}{2})}
  \bigl(w_j(s)v_i(s) - v_j(s)w_i(s)\bigr)\d s\ ,
\end{equation}
where $i = 1,2,3$. The main result of this section states that, if the 
initial data are small, \reff{SV3int} has global solutions in $\L^2(m)$ 
which decay exponentially to zero as $\tau \to +\infty$.

\begin{theorem}\label{zeroasym} Let $0 < \mu \le 1$ and $m > 2\mu+
\frac{1}{2}$. There exists $r_0 > 0$ such that, for all initial data 
$\ww_0 \in \L^2(m)$ with $\|\ww_0\|_m \le r_0$, equation \reff{SV3int} 
has a unique global solution $\ww \in C^0([0,\infty),\L^2(m))$ 
satisfying $\ww(0) = \ww_0$. In addition, there exists $K_0 \ge 1$ 
such that
\begin{equation}\label{wdecay}
  \|\ww(\tau)\|_m \le K_0 \,\e^{-\mu\tau}\|\ww_0\|_m\ , \quad \tau \ge 0\ .
\end{equation}
\end{theorem}

\proof Given $\ww_0 \in \L^2(m)$, we shall solve \reff{SV3int} in the 
Banach space
$$
  X = \bigl\{\ww \in C^0([0,+\infty),\L^2(m)) \,\big|\, 
  \|\ww\|_X = \sup_{\tau\ge 0}\|\ww(\tau)\|_m \e^{\mu\tau} < \infty
  \bigr\}\ .
$$
We first note that $\tau \mapsto \e^{\tau\Lambda}\ww_0 \in X$; 
namely, there exists $C_1 \ge 1$ such that
\begin{equation}\label{linbd}
  \|\e^{\tau\Lambda}\ww_0\|_m \le C_1 \e^{-\mu\tau}\|\ww_0\|_m\ ,
  \quad \tau \ge 0\ .
\end{equation}
Indeed, if $m \le \frac{5}{2}$, \reff{linbd} is nothing but \reff{Sg1}
with $\epsilon = m-2\mu-\frac{1}{2}$, $\alpha = 0$, and $n = -1$ or $0$.  
If $m > \frac{5}{2}$, \reff{linbd} follows from \reff{Sg2} with 
$\alpha = 0$ and $n = 0$. 

Next, given $\ww \in C^0([0,+\infty),\L^2(m))$, we define $\FF[\ww]
\in C^0([0,+\infty),\L^2(m))$ by
\begin{equation}\label{FFdef}
  (F_i[\ww])(\tau) = \sum_{j=1}^3 \int_0^\tau \partial_j 
  \e^{(\tau-s)(\Lambda-\frac{1}{2})}\bigl(w_j(s)v_i(s) - 
  v_j(s)w_i(s)\bigr)\d s\ , \quad \tau \ge 0\ .
\end{equation}
We shall prove that $\FF$ maps $X$ into $X$, and that there exists
$C_2 > 0$ such that
\begin{equation}\label{nonbd}
  \|\FF[\ww]\|_X \le C_2 \|\ww\|_X^2\ ,\quad 
  \|\FF[\ww]-\FF[\tilde \ww]\|_X \le C_2 \|\ww{-}\tilde \ww\|_X
  (\|\ww\|_X + \|\tilde \ww\|_X)\ ,
\end{equation}
for all $\ww, \tilde \ww \in X$. As is easily verified, the bounds
\reff{linbd}, \reff{nonbd} imply that the map $\ww \mapsto 
\e^{\tau\Lambda}\ww_0 + \FF[\ww]$ has a unique fixed point in 
the ball $\{\ww \in X\,|\,\|\ww\|_X \le R\}$ if $R < (2C_2)^{-1}$
and $\|\ww_0\|_m \le (2C_1)^{-1}R$. Using Gronwall's lemma, it is
then straightforward to show that this fixed point is actually the
unique solution of \reff{SV3int} in the space $C^0([0,+\infty),\L^2(m))$. 
Finally, since $\|\ww\|_X \le C_1 \|\ww_0\|_m + C_2\|\ww\|_X^2 \le
C_1 \|\ww_0\|_m + \frac{1}{2}\|\ww\|_X$, the bound \reff{wdecay}
holds with $K_0 = 2C_1$. 

To prove \reff{nonbd}, we use the following estimate, which is a 
consequence of Propositions~\ref{Sgestim} and \ref{LpLqestim}. Assume that
$f : \real^3 \to \real$ satisfies $\rho^m f \in L^{3/2}(\real^3)$, 
where $\rho(\xi) = 1+|\xi|$. Then there exists $C_3 > 0$ such that, 
for $j = 1,2,3$,
\begin{equation}\label{m32}
  \|\partial_j \e^{\tau\Lambda}f\|_m \le C_3 
  \frac{\e^{-\nu\tau}}{a(\tau)^{3/4}}\,|\rho^m f|_{3/2}\ , 
  \quad \tau > 0\ ,
\end{equation}
where $a(\tau) = 1-\e^{-\tau}$ and $\nu = \min(\mu,\frac{1}{2})$. 
Indeed, if $0 < \tau < 2$, then \reff{m32} follows from 
\reff{LpLq} with $p = 2$ and $q = 3/2$. If $\tau \ge 2$, we have
\begin{equation}\label{auxin}
  \|\partial_j \e^{(\tau-1)\Lambda} \e^{\Lambda}f\|_m \le
  C \e^{-\nu(\tau-1)}\|\e^\Lambda f\|_m \le C' \e^{-\nu\tau}
  |\rho^m f|_{3/2}\ ,
\end{equation}
where the second inequality is again a consequence of \reff{LpLq}. 
The first inequality in \reff{auxin} follows from \reff{Sg1} with 
$\epsilon = m-2\mu-\frac{1}{2}$ and $n = -1$ if $m \le 3/2$, and  
from \reff{Sg2} with $n = -1$ if $m > 3/2$. This proves \reff{m32}
for all $\tau > 0$. 

Given $\ww \in X$ and $s \ge 0$, we apply \reff{m32} to $f = 
w_j(s)v_i(s)-v_j(s)w_i(s)$, where $i,j \in \{1,2,3\}$ and 
$\vv = (v_1,v_2,v_3)$ is the velocity field obtained from $\ww$
via the Biot-Savart law. Using H\"older's inequality and Lemma~\ref{HLS3}, 
we can bound
\begin{equation}\label{m33}
  |\rho^m w_j v_i|_{3/2} \le |\rho^m w_j|_2 |v_i|_6 
  \le C \|w_j\|_m |w_i|_2 \le C\|w_j\|_m \|w_i\|_m\ ,
\end{equation}
so that $|\rho^m f|_{3/2} \le C_4 \|\ww(s)\|_m^2$ for some $C_4 > 0$. 
 Combining this bound with \reff{m32} and using \reff{FFdef}, we obtain,
for all $\tau > 0$, 
\begin{eqnarray*}
  \|F_i([\ww])(\tau)\|_m &\le& 3C_3 C_4 \int_0^\tau 
    \frac{\e^{-(\nu+\frac{1}{2})(\tau-s)}}{a(\tau{-}s)^{3/4}}
    \|\ww(s)\|_m^2 \d s\\
  &\le& 3C_3 C_4 \|\ww\|_X^2 \int_0^\tau 
    \frac{\e^{-(\nu+\frac{1}{2})(\tau-s)}}{a(\tau{-}s)^{3/4}}
    \,\e^{-2\mu s}\d s \le C_2 \|\ww\|_X^2 \e^{-\mu\tau}\ ,
\end{eqnarray*}
since $\nu+\frac{1}{2} \ge \mu$. This establishes the first inequality 
in \reff{nonbd}, and the second one can be proved along the same lines. 
The proof of Theorem~\ref{zeroasym} is now complete. \QED

\medskip Since the semigroup $\e^{\tau\Lambda}$ is not analytic in
$\L^2(m)$, the solution $\ww$ given by Theorem~\ref{zeroasym}
is in general not a smooth function of $\tau$. In particular,
$\tau \mapsto \ww(\tau) \notin C^1((0,+\infty),\L^2(m))$, so that
$\ww$ is not a classical solution of \reff{SV3} in $\L^2(m)$. 
Nevertheless, following the common use, we shall often refer to 
$\ww$ as to the (mild) solution of \reff{SV3} in $\L^2(m)$. 
Remark that the evolution defined by \reff{SV3} is regularizing
in the sense that $\ww(\xi,\tau)$ is a smooth function of $\xi \in 
\real^3$ for any $\tau > 0$. This property is well-known, and will
not be proved here. We only quote the following result:

\begin{proposition}\label{wregul} Let $0 < \mu \le 1$, $m > 2\mu+
\frac{1}{2}$, and let $\ww \in C^0([0,\infty),\L^2(m))$ be the solution 
of \reff{SV3} given by Theorem~\ref{zeroasym}.
There exists $K_1 > 0$ such that, for all $p \in [2,+\infty]$, 
\begin{equation}\label{wlpdecay}
  |\rho^m \ww(\tau)|_p \le K_1 (1+\tau^{-\gamma_p})\,\e^{-\mu\tau}
  \|\ww_0\|_m\ , \quad \tau > 0\ ,
\end{equation}
where $\rho(\xi) = 1+|\xi|$ and $\gamma_p = \frac{3}{2}(\frac{1}{2}-
\frac{1}{p})$. 
\end{proposition}

\proof In view of \reff{wdecay}, it is clearly sufficient to prove
\reff{wlpdecay} for $0 < \tau \le 1$. This can be done by a standard
bootstrap argument, using Proposition~\ref{LpLqestim}, 
Lemma~\ref{HLS3}, and the integral equation \reff{SV3int} satisfied 
by $\ww$. We omit the details. \QED

\begin{corollary}\label{zerolp}  Let $0 < \mu \le 1$ and $m > 2\mu+
\frac{1}{2}$. Let $\ww \in C^0([0,\infty),\L^2(m))$ be the solution of 
\reff{SV3} given by Theorem~\ref{zeroasym}, and let $\vv$ be the 
corresponding velocity field. There exists $K_2 > 0$ such that, for all 
$\tau \ge 1$, 
\begin{equation}\label{estA}
   |\ww(\tau)|_p \le K_2 \e^{-\mu\tau}\|\ww_0\|_m \quad
   \hbox{for all } \cases{p \in (p_m,\infty] & if~ 
   $\frac{1}{2} < m \le \frac{3}{2}$\ ,\cr p \in [1,\infty] & 
   if~ $m > \frac{3}{2}$\ ,}
\end{equation}
and
\begin{equation}\label{estB}
   |\vv(\tau)|_q \le K_2 \e^{-\mu\tau}\|\ww_0\|_m \quad
   \hbox{for all } \cases{q \in (q_m,\infty] & if~ 
   $\frac{1}{2} < m \le \frac{5}{2}$\ ,\cr q \in (1,\infty] & 
   if~ $m > \frac{5}{2}$\ ,}
\end{equation}
where $p_m = \frac{6}{3+2m}$ and $q_m = \frac{6}{1+2m}$. Moreover, 
the bounds \reff{estA}, \reff{estB} hold for all $\tau \ge 0$ if $p \le 2$, 
$q \le 6$.
\end{corollary}

\proof If $p > 2$, \reff{estA} is a direct consequence of \reff{wlpdecay}. 
If $p_m < p \le 2$, then $|\ww|_p \le C|\rho^m \ww|_2$ by H\"older's 
inequality, and \reff{estA} again follows from \reff{wlpdecay}. 
Using \reff{estA} and Lemma~\ref{HLS3}, we easily obtain \reff{estB}
if $q > 3/2$. Finally, if $m > 3/2$ and $\max(1,q_m) < q \le 3/2$, 
we can pick $n < \min(m,5/2)$ such that $q > q_n$. Using H\"older's 
inequality and Proposition~\ref{velvort3}, we find $|\vv|_q \le 
C|\rho^n\vv|_6 \le C|\rho^n \ww|_2 \le C\|\ww\|_m$, and \reff{estB}
follows. \QED

\begin{remark} The fact that the value $q = 1$ is excluded in \reff{estB} 
if $m > \frac{5}{2}$ is not a technical restriction. As is shown in 
Corollary~\ref{3velint}, the velocity field $\vv(\xi,\tau)$ is not integrable 
in this case, unless $\int_{\real^3}\xi_i w_j(\xi,\tau)\d\xi = 0$ for all 
$i,j \in \{1,2,3\}$. 
\end{remark}

For the vorticity $\oomega(x,t)$ and the velocity $\uu(x,t)$ in 
the original variables, Corollary~\ref{zerolp} implies, for the 
same values of $p$ and $q$:
$$
   |\oomega(t)|_p \le K_2 t^{-1-\mu+\frac{3}{2p}}\|\oomega_0\|_m\ ,
   \quad |\uu(t)|_q \le K_2 t^{-\frac{1}{2}-\mu+\frac{3}{2q}}
   \|\oomega_0\|_m\ , \quad t \ge 1\ .
$$

\medskip
We now explain our motivation for introducing the scaling variables 
\reff{scalvar}. As is shown in Appendix~\ref{spectrum}, 
the spectrum of $\Lambda$ acting on $\L^2(m)$ can be decomposed as
$\sigma(\Lambda) = \sigma_d(\Lambda) \cup \sigma_c(\Lambda)$, where
$$
  \sigma_d(\Lambda) = \Bigl\{-\frac{k+1}{2} \,\Big|\, 
  k \in \intplus^*\Bigr\}\ , \quad 
  \sigma_c(\Lambda) = \Bigl\{\lambda \in \complex\,\Big|\, 
  \Re(\lambda) \le \frac{1}{4} - \frac{m}{2}\Bigr\}\ .
$$
(See Fig.~1.)
Remark that the discrete spectrum $\sigma_d(\Lambda)$ does not depend on
$m$, whereas the continuous spectrum $\sigma_c(\Lambda)$ can be shifted
arbitrarily far away from the origin by choosing $m$ appropriately. 
Therefore, if $m \ge 0$ is sufficiently large, the long-time 
behavior of the solutions of \reff{SV3} in a neighborhood of the
origin is governed by a {\it finite} system of ordinary differential
equations. This system is obtained by projecting \reff{SV3} onto the 
finite-dimensional subspace of $\L^2(m)$ spanned by the eigenfunctions 
of $\Lambda$ corresponding to the first eigenvalues $\lambda_k = 
-\frac{k+1}{2}$, with $k = 1,2,\dots,k_0$. 

\figurewithtex 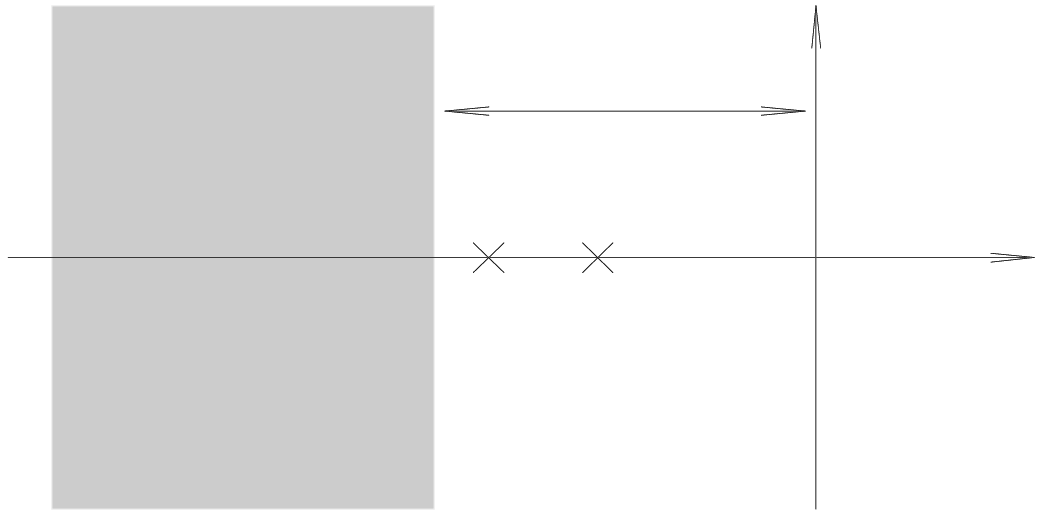 Fig1.tex 5.5 11 {\bf Fig.~1:} The spectrum
of the linear operator $\Lambda$ in $\L^2(m)$, when $m = 4$.\cr

A rigorous justification of this reduction, using invariant manifold
theory, can be found in \cite{gallay:2001} for the two-dimensional 
vorticity equation. Specifically, given any $\nu > 0$, we prove the 
existence of a finite-dimensional locally invariant manifold, which
is tangent at the origin to the spectral subspace corresponding to the
first eigenvalues of $\Lambda$, and which is approached at a rate 
$\cO(\e^{-\nu\tau})$ or faster by any solution of \reff{SV3} that stays 
in a neighborhood of the origin for all times. This method allows, at 
least in principle, to compute the long-time asymptotics of the solutions 
to arbitrarily high order by studying a finite-dimensional dynamical 
system -- the restriction of the rescaled vorticity equation to the manifold. 

Invariant manifolds can be constructed in the three-dimensional
case also, and we use them in our discussion of the set of
solutions of the Navier-Stokes equations which decay
``faster than expected'' in Section \ref{ssmanif}.  However,
for computing the asymptotics, we show that one can also use
a different approach.  Given $k_0 \in \intplus^*$ and 
$m > k_0 + \frac{3}{2}$, we decompose any solution $\ww$ of \reff{SV3}
in $\L^2(m)$ as
$$
   \ww(\xi,\tau) = \sum_{k=1}^{k_0} \sum_{\ell=1}^{k(k{+}2)} 
   \alpha_{k\ell}(\tau)\ww_{k\ell}(\xi) + \RR(\xi,\tau)\ ,
$$
where $\alpha_{k\ell} \in \real$ and, for any $k \in \{1,\dots,k_0\}$, 
$\{\ww_{k\ell}\,|\, \ell = 1,\dots,k(k{+}2)\}$ is a basis of the eigenspace 
$\{\ww \in \L^2(m) \,|\, \Lambda\ww =-\frac{k+1}{2}\ww\}$. Using this 
decomposition, \reff{SV3} becomes a system of ordinary differential 
equations for the coefficients $\alpha_{k\ell}$ coupled to a partial 
differential equation for the remainder $\RR$. A direct analysis of this 
system allows to compute the asymptotics up to order $\cO(\e^{-\nu\tau})$, 
where $\nu = \min(\frac{k_0+2}{2},\frac{m}{2}-\frac{1}{4})$. 
This program is carried out in Section~\ref{first} for $k_0 = 1$ 
(first order asymptotics) and in Section~\ref{second} for $k_0 = 2$ 
(second order asymptotics). 

%%%%%%%%%%%%%%%%%%%%%%%%%%%%%%%%%%%%%%%%%%%%%%%%%%%%%%%%%%%%%%%%%%%%%%%%%%%

\section{First-order asymptotics}\label{first}

In this section, we consider the behavior of the solutions of \reff{SV3} 
in $\L^2(m)$ with $\frac{5}{2} < m \le \frac{7}{2}$. In this space, 
the discrete spectrum of $\Lambda$ consists of a single isolated 
eigenvalue $\lambda_1 = -1$, of multiplicity $3$. A convenient basis
of eigenvectors is given by $\{\ff_1,\ff_2,\ff_3\}$, where $\ff_i = 
\rot(G\ee_i)$. Here and in the sequel, $G$ is the Gaussian function
$$
   G(\xi) = \frac{1}{(4\pi)^{3/2}}\,\e^{-|\xi|^2/4}\ ,
   \quad \xi \in \real^3\ ,
$$
and $\{\ee_1,\ee_2,\ee_3\}$ denotes the canonical basis of $\real^3$.
A short calculation shows that $\ff_i = \pp_i G$ for $i=1,2,3$, where
$\pp_i(\xi) = \frac{1}{2}(\ee_i\wedge \xi)$. Explicitly, 
\begin{equation}\label{ppdef}
  \pp_1(\xi) = \frac{1}{2}\pmatrix{0 \cr -\xi_3 \cr \xi_2}\ , \quad
  \pp_2(\xi) = \frac{1}{2}\pmatrix{\xi_3 \cr 0 \cr -\xi_1}\ , \quad
  \pp_3(\xi) = \frac{1}{2}\pmatrix{-\xi_2 \cr \xi_1 \cr 0}\ .
\end{equation}
The vector fields $\pp_i$ satisfy $\div \pp_i = 0$ and $\rot 
\pp_i = \ee_i$. Integrating by parts, we thus find
$$
  \int_{\real^3}\pp_i \cdot \ff_j \d\xi = \int_{\real^3}\rot(\pp_i)
  \cdot(G\ee_j)\d\xi = (\ee_i\cdot\ee_j)\int_{\real^3}G \d\xi  = 
  \delta_{ij}\ .
$$
Moreover, if $\Lambda^* = \Delta -\frac{1}{2}(\xi\cdot\nabla)-\frac{1}{2}$ 
is the formal adjoint of $\Lambda$, it is easy to verify that 
$\Lambda^* \pp_i = -\pp_i$ for $i = 1,2,3$. The velocity fields 
$\vv^{\ff_i}$ corresponding to $\ff_i$ are computed in 
Appendix~\ref{velocity3d}. In particular, we mention that 
$|\vv^{\ff_i}(\xi)| \sim |\xi|^{-3}$ as $|\xi| \to \infty$, so that 
$\vv^{\ff_i} \in \L^q(\real^3)$ for all $q > 1$. Using these notations, 
any solution $\ww$ of \reff{SV3} in $\L^2(m)$ can be decomposed as
\begin{equation}\label{wdecomp}
  \ww(\xi,\tau) = \sum_{i=1}^3 \beta_i(\tau)\ff_i(\xi) + \RR(\xi,\tau)\ ,
\end{equation}
where
\begin{equation}\label{betadef}
  \beta_i(\tau) = \int_{\real^3} \pp_i(\xi)\cdot\ww(\xi,\tau)\d\xi\ ,
  \quad i = 1,2,3\ .
\end{equation}
Then $\RR(\cdot,\tau)$ belongs to the subspace $\cW_1$ of $\L^2(m)$ 
defined in \reff{Wndef}, which is also the spectral subspace 
associated with the continuous spectrum $\{\lambda \in \complex\,|\,
\Re(\lambda) \le \frac{1}{4}-\frac{m}{2}\}$ of the operator $\Lambda$. 
As in the two-dimensional case \cite{gallay:2001}, the coefficients 
$\beta_i$ obey a linear evolution equation:

\begin{lemma}\label{3betadot}
Assume that $m > \frac{5}{2}$, and let $\ww \in C^0([0,T],\L^2(m))$ be
a solution of \reff{SV3}. Then the coefficients $\beta_i$ defined by 
\reff{betadef} satisfy, for all $\tau \in [0,T]$,
$$
   \dot \beta_i(\tau) = -\beta_i(\tau)\ ,\quad i = 1,2,3\ .
$$
\end{lemma}

\proof Since $\rot(\vv\wedge\ww) = (\ww\cdot\nabla)\vv - (\vv\cdot\nabla)
\ww$, \reff{SV3} is equivalent to
$$
  \partial_\tau \ww = \Lambda \ww + \rot(\vv\wedge\ww)\ , \quad \div \ww 
  = 0\ .
$$
Differentiating \reff{betadef} formally with respect to $\tau$ and 
integrating by parts, we thus find
\begin{equation}\label{auxdot}
  \dot \beta_i = \int_{\real^3} \pp_i\cdot\bigl(\Lambda \ww
  + \rot(\vv\wedge\ww)\bigr)\d\xi 
  = -\beta_i + \int_{\real^3} \ee_i \cdot (\vv\wedge\ww)\d\xi\ ,
\end{equation}
where we used the fact that $\Lambda^* \pp_i = -\pp_i$. Since the
right-hand side of \reff{auxdot} belongs to $C^0([0,T])$ and depends
continuously on $\ww$, the calculations above can be justified by a 
density argument. In particular, $\beta_i \in C^1([0,T])$ for 
$i = 1,2,3$. Finally, using the identity $\vv\wedge\ww = \vv\wedge\rot\vv = 
\frac{1}{2}\nabla |\vv|^2 - (\vv\cdot\nabla)\vv$ and the fact that 
$\div\vv = 0$, we see that the last integral in \reff{auxdot} vanishes, 
hence $\dot\beta_i = -\beta_i$. \QED

\medskip
In particular, it follows from Lemma~\ref{3betadot} that the 
subspace $\cW_1$ is {\it invariant} under the evolution defined 
by \reff{SV3}. The remainder $\RR$ in \reff{wdecomp} satisfies the 
equation
$$
  \partial_\tau \RR = \Lambda \RR + Q_1\bigl((\ww\cdot\nabla)\vv - 
  (\vv\cdot\nabla)\ww\bigr)\ , \quad \div \RR = 0\ ,
$$
where $Q_1 : \L^2(m) \to \cW_1$ is the spectral projection (for the
operator $\Lambda$) onto the subspace $\cW_1$. Explicitly, 
$$
  Q_1 \ww = \ww - \sum_{i=1}^3 \Bigl(\int_{\real^3}\pp_i\cdot
  \ww \d\xi\Bigr) \ff_i\ .
$$

The following result describes the first order asymptotics 
of $\ww(\tau)$ as $\tau \to \infty$. 

\begin{theorem}\label{firstasym}
Let $1 < \nu \le \frac{3}{2}$, $m > 2\nu+\frac{1}{2}$, and let 
$\ww \in C^0([0,\infty),\L^2(m))$ be the solution of \reff{SV3}
given by Theorem~\ref{zeroasym} with $\mu = 1$. Then there exists 
$K_3 \ge 1$ such that
$$
   \bigl\|\ww(\tau) - \sum_{i=1}^3 b_i \e^{-\tau}\ff_i\bigr\|_m 
   \le K_3 \e^{-\nu \tau} \|\ww_0\|_m\ , \quad \tau \ge 0\ ,
$$
where $b_i = \int_{\real^3}\pp_i\cdot\ww_0\d\xi$, $i = 1,2,3$. 
\end{theorem}

\proof If $\ww \in C^0([0,\infty),\L^2(m))$ is the solution of \reff{SV3}
given by Theorem~\ref{zeroasym} and $\vv$ the corresponding velocity
field, we define $\beta_i$ and $\RR$ by \reff{wdecomp}, \reff{betadef}.
By Lemma~\ref{3betadot}, $\beta_i(\tau) = b_i \e^{-\tau}$ for $i = 1,2,3$.
To bound the remainder $\RR$, we use the integral equation
\begin{equation}\label{Rint}
   \RR(\tau) = \e^{\tau\Lambda}\RR_0 + \int_0^\tau Q_1 
   \PPhi(\tau-s,\ww(s),\vv(s))\d s\ , 
\end{equation}
where $\RR_0 = Q_1 \ww_0$ and $\PPhi(\sigma,\ww,\vv) = \e^{\sigma\Lambda}
((\ww\cdot\nabla)\vv - (\vv\cdot\nabla)\ww)$. 

By Proposition~\ref{Sgestim}, there exists $C_1 > 0$ such that
\begin{equation}\label{linr}
  \|\e^{\tau\Lambda}\RR_0\|_m \le C_1 \e^{-\nu\tau}\|\ww_0\|_m\ , 
  \quad \tau \ge 0\ .
\end{equation}
Indeed, if $\frac{5}{2} < m \le \frac{7}{2}$, \reff{linr} is a 
consequence of \reff{Sg1} with $\alpha = 0$, $n = 1$, and 
$\epsilon = m-2\nu-\frac{1}{2}$. If $m > \frac{7}{2}$, \reff{linr}
follows from \reff{Sg2} with $\alpha = 0$ and $n = 1$. 

To estimate the integral in \reff{Rint}, we proceed as in the proof of
Theorem~\ref{zeroasym}. Exchanging $\nabla$ with $\e^{\sigma\Lambda}$, 
we can write $\PPhi = (\Phi_1,\Phi_2,\Phi_3)$, where
$$
  \Phi_i(\sigma,\ww,\vv) = \sum_{j=1}^3 \partial_j \e^{\sigma
  (\Lambda-\frac{1}{2})}(w_j v_i - v_j w_i)\ , \quad i = 1,2,3\ ,
$$
see \reff{SV3int}. By \reff{m32}, \reff{m33}, there exists $C_2 > 0$ such 
that $\|Q_1\PPhi(\sigma,\ww,\vv)\|_m \le C_2 \sigma^{-3/4}\|\ww\|_m^2$ for 
all $\sigma \in [0,1]$. Using \reff{linr} and \reff{wdecay}, we thus
find
$$
   \|\RR(\tau)\|_m \le C_1 \|\ww_0\|_m + \int_0^\tau C_2(\tau{-}s)^{-3/4}
   (K_0 \|\ww_0\|_m)^2 \d s \le C \|\ww_0\|_m\ , 
$$
for all $\tau \in [0,1]$. (Here and in the sequel, we use the fact 
that $\|\ww_0\|_m \le r_0$, see Theorem~\ref{zeroasym}.) If $\tau > 1$, 
we write $\RR(\tau) = \e^{\tau\Lambda}\RR_0 + \RR_1(\tau) +\RR_2(\tau)$, 
where
\begin{eqnarray*}
  \RR_1(\tau) &=& \int_0^{\tau-1}\e^{(\tau-1-s)\Lambda}Q_1 
  \PPhi(1,\ww(s),\vv(s))\d s\ , \\
  \RR_2(\tau) &=& \int_{\tau-1}^\tau Q_1 \PPhi(\tau-s,\ww(s),\vv(s))
  \d s\ .
\end{eqnarray*}
Proceeding as above, we obtain
\begin{eqnarray*}
  \|\RR_1(\tau)\|_m &\le& \int_0^{\tau-1} C_1 \e^{-\nu(\tau-1-s)}
   C_2 (K_0 \e^{-s}\|\ww_0\|_m)^2 \d s \le C \e^{-\nu \tau}\|\ww_0\|_m^2\ ,\\
 \|\RR_2(\tau)\|_m &\le& \int_{\tau-1}^\tau C_2 (\tau{-}s)^{-3/4}
  (K_0 \e^{-s}\|\ww_0\|_m)^2 \d s \le C \e^{-2\tau}\|\ww_0\|_m^2\ .
\end{eqnarray*}
Thus, there exists $K_3 > 0$ such that $\|\RR(\tau)\|_m \le K_3 
\e^{-\nu\tau}\|\ww_0\|_m$, for all $\tau \ge 0$. \QED

\begin{corollary}\label{firstlp} Let $1 < \nu \le \frac{3}{2}$ and
$m > 2\nu+\frac{1}{2}$. Let $\ww \in C^0([0,\infty),\L^2(m))$ 
be the solution of \reff{SV3} given by Theorem~\ref{zeroasym} with
$\mu = 1$, and let $\vv$ be the corresponding velocity field. There exists 
$K_4 > 0$ such that, for all $\tau \ge 1$, 
\begin{eqnarray}\label{esttA}
  &&\Big|\ww(\tau) - \sum_{i=1}^3 b_i \e^{-\tau}\ff_i\Big|_p \le K_4 
  \e^{-\nu\tau}\|\ww_0\|_m\ , \quad 1 \le p \le \infty\ ,\\ \label{esttB}
  &&\Big|\vv(\tau) - \sum_{i=1}^3 b_i \e^{-\tau}\vv^{\ff_i}\Big|_q \le K_4 
  \e^{-\nu\tau}\|\ww_0\|_m\ , \quad 1 \le q \le \infty\ ,
\end{eqnarray}
where $\vv^{\ff_i}$ is given by \reff{vfidef}. Moreover, the bounds 
\reff{esttA}, \reff{esttB} hold for all $\tau \ge 0$ if $p \le 2$, $q \le 6$.
\end{corollary}

\proof Using the analogue of Proposition~\ref{wregul} for $\RR(\xi,\tau)$
and proceeding as in the proof of Corollary~\ref{zerolp}, we obtain 
\reff{esttA} for $1 \le p \le \infty$ and \reff{esttB} for $1 < q \le 
\infty$. To prove \reff{esttB} for $q = 1$, we pick $n \in \real$ 
such that $5/2 < n < \min(m,7/2)$. If $\vv^{\RR}$ is the velocity field 
obtained from $\RR$ via the Biot-Savart law \reff{SBS3}, then using 
H\"older's inequality and Proposition~\ref{velvort3} we can bound 
$|\vv^{\RR}|_1 \le C|\rho^n \vv^{\RR}|_6 \le C|\rho^n \RR|_2 \le C\|\RR\|_m$, 
and the result follows. \QED

\medskip
In terms of the original variables, Corollary~\ref{firstlp} shows that, 
for all $t \ge 1$, 
\begin{eqnarray*}
  &&|\oomega(t) - \oomega_{app}(t)|_p \le C t^{-1-\nu+\frac{3}{2p}}
  \|\oomega_0\|_m\ , \quad 1 \le p \le \infty\ ,\\
  &&|\uu(t) - \uu_{app}(t)|_q \le C t^{-\frac{1}{2}-\nu+\frac{3}{2q}}
  \|\oomega_0\|_m\ , \quad 1 \le q \le \infty\ ,
\end{eqnarray*}
where $\oomega_{app}(x,t)$, $\uu_{app}(x,t)$ are the self-similar 
vector fields defined by
$$
   \oomega_{app}(x,t) = \sum_{i=1}^3 \frac{b_i}{(1+t)^2}
   \,\ff_i\Bigl(\frac{x}{\sqrt{1+t}}\Bigr)\ , \quad 
   \uu_{app}(x,t) = \sum_{i=1}^3 \frac{b_i}{(1+t)^{3/2}}
   \,\vv^{\ff_i}\Bigl(\frac{x}{\sqrt{1+t}}\Bigr)\ .
$$
Remark that $\uu_{app}(\cdot,t) \notin \L^1(\real^3)$, unless
$b_1 = b_2 = b_3 = 0 $. 

%%%%%%%%%%%%%%%%%%%%%%%%%%%%%%%%%%%%%%%%%%%%%%%%%%%%%%%%%%%%%%%%%%%%%%%%%%%

\section{Second-order asymptotics}\label{second}

We now turn our attention to the solutions of \reff{SV3} in $\L^2(m)$
with $\frac{7}{2} < m \le \frac{9}{2}$. Acting on this space, the operator
$\Lambda$ has exactly two isolated eigenvalues: $\lambda_1 = -1$
(of multiplicity $3$) and $\lambda_2 = -\frac{3}{2}$ (of multiplicity $8$). 
Let $\cE_2$ be the subspace of $\L^2(m)$ spanned by the eigenfunctions 
corresponding to $\lambda_2$, see Appendix~\ref{spectrum}. A convenient 
basis of $\cE_2$ is provided by the vector fields $\gg_i$ and $\hh_{ij}$, 
which we now define. 

\smallskip{\noindent}
{\bf a)} For $i=1,2,3$, let $\gg_i = \rot\ff_i = \rot(\pp_i G) = 
\frac{G}{4}((4-|\xi|^2)\ee_i + \xi_i\xi)$. Explicitly, 
\begin{equation}\label{ggdef}
  \gg_1 = \frac{G}{4}
    \pmatrix{4{-}\xi_2^2 {-}\xi_3^3 \cr \xi_1\xi_2 \cr \xi_1\xi_3}\ , \quad
  \gg_2 = \frac{G}{4}
    \pmatrix{\xi_1\xi_2 \cr 4{-}\xi_1^2{-}\xi_3^2 \cr \xi_2\xi_3}\ , \quad
  \gg_3 = \frac{G}{4}
    \pmatrix{\xi_1\xi_3 \cr \xi_2\xi_3 \cr 4{-}\xi_1^2{-}\xi_2^2}\ .
\end{equation}
Then $\div \gg_i = 0$ and $\Lambda \gg_i = -\frac{3}{2}\gg_i$. 
By construction, the velocity field associated to $\gg_i$ is $\vv^{\gg_i} 
\equiv \ff_i$. In particular, $\vv^{\gg_i}$ has a Gaussian decay as 
$|\xi| \to \infty$. We also define
\begin{equation}\label{qqdef}
  \qq_1(\xi) = \frac{1}{2}
    \pmatrix{2-\xi_1^2 \cr \xi_1\xi_2 \cr \xi_1\xi_3}\ , \quad
  \qq_2(\xi) = \frac{1}{2}
    \pmatrix{\xi_1\xi_2 \cr 2-\xi_2^2 \cr \xi_2\xi_3}\ , \quad
  \qq_3(\xi) = \frac{1}{2}
    \pmatrix{\xi_1\xi_3 \cr \xi_2\xi_3 \cr 2-\xi_3^2}\ .
\end{equation}
Then $\div\qq_i = 0$, $\rot \qq_i = \pp_i$, and $\Lambda^* \qq_i = 
-\frac{3}{2}\qq_i$. It follows that
$$
   \int_{\real^3} \qq_i \cdot\gg_j \d\xi = \int_{\real^3} \rot(\qq_i)\cdot
   \ff_j \d\xi = \int_{\real^3} \pp_i \cdot \ff_j \d\xi = \delta_{ij}\ .
$$

\smallskip{\noindent}
{\bf b)} For $(ij) \in S = \{(11),(12),(13),(22),(23)\}$, we define
$\hh_{ij} = \partial_i \ff_j + \partial_j \ff_i$. Explicitly, we
have $\hh_{ii} = -\xi_i\ff_i$ for $i = 1,2$ and
\begin{equation}\label{hhdef}
  \hh_{12} = \frac{G}{4}
    \pmatrix{-\xi_1\xi_3 \cr \xi_2\xi_3 \cr \xi_1^2-\xi_2^2}\ ,\quad
  \hh_{13} = \frac{G}{4}
    \pmatrix{\xi_1\xi_2 \cr \xi_3^2-\xi_1^2 \cr -\xi_2\xi_3}\ ,\quad
  \hh_{23} = \frac{G}{4}
    \pmatrix{\xi_2^2-\xi_3^2 \cr -\xi_1\xi_2 \cr \xi_1\xi_3}\ .
\end{equation}
Then $\div \hh_{ij} = 0$ and $\Lambda \hh_{ij} = -\frac{3}{2}\hh_{ij}$ 
for all $(ij) \in S$. The velocity fields $\vv^{\hh_{ij}}$ corresponding
to $\hh_{ij}$ are computed in Appendix~\ref{velocity3d}. In particular, 
we remark that $|\vv^{\hh_{ij}}(\xi)| \sim |\xi|^{-4}$ as $|\xi| \to \infty$, 
so that $\rho \vv^{\hh_{ij}} \notin L^1(\real^3)^3$. We also define 
$\rr_{11} = \frac{1}{2}\xi_1\xi_3\ee_2$, $\rr_{22} = -\frac{1}{2}
\xi_2\xi_3\ee_1$, and 
\begin{equation}\label{rrdef}
  \rr_{12} = \frac{1}{2}
    \pmatrix{-\xi_1\xi_3 \cr \xi_2\xi_3 \cr 0}\ ,\quad
  \rr_{13} = \frac{1}{2}
    \pmatrix{\xi_1\xi_2 \cr 0 \cr -\xi_2\xi_3}\ ,\quad
  \rr_{23} = \frac{1}{2}
    \pmatrix{0 \cr -\xi_1\xi_2 \cr \xi_1\xi_3}\ .
\end{equation}
Then $\div \rr_{ij} = 0$ and $\Lambda^* \rr_{ij} = -\frac{3}{2}\rr_{ij}$ 
for all $(ij) \in S$. A direct calculation shows that the following 
orthogonality relations are satisfied:
$$
   \int_{\real^3} \rr_{ij}\cdot\hh_{kl}\d\xi = \delta_{ik}\delta_{jl}\ ,
   \quad \int_{\real^3} \rr_{ij}\cdot\gg_k \d\xi = \int_{\real^3} 
   \qq_k \cdot \hh_{ij}\d\xi = 0\ .
$$

Using these notations, any solution $\ww$ of \reff{SV3} in $\L^2(m)$ 
can be decomposed as
\begin{equation}\label{decompp}
  \ww(\xi,\tau) = \sum_{i=1}^3 \beta_i(\tau)\ff_i(\xi) + 
  \sum_{i=1}^3 \gamma_i(\tau)\gg_i(\xi) + \sum_{(ij)\in S} 
  \zeta_{ij}(\tau)\hh_{ij}(\xi) + \RR(\xi,\tau)\ ,
\end{equation}
where $\beta_i(\tau)$ is given by \reff{betadef} and 
\begin{equation}\label{gamzetdef}
  \gamma_i(\tau) = \int_{\real^3} \qq_i(\xi)\cdot\ww(\xi,\tau)\d\xi\ ,
  \quad \zeta_{ij}(\tau) = \int_{\real^3} \rr_{ij}(\xi)\cdot
  \ww(\xi,\tau)\d\xi\ .
\end{equation}
Then $\RR(\cdot,\tau)$ belongs to the subspace $\cW_2$ of $\L^2(m)$
defined in \reff{Wndef}, which coincides with the spectral subspace
associated with the continuous spectrum $\{\lambda \in \complex\,|\,
\Re(\lambda) \le \frac{1}{4}-\frac{m}{2}\}$ of the operator $\Lambda$. 
The coefficients $\gamma_i$ and $\zeta_{ij}$ satisfy the following 
evolution equations:

\begin{lemma}\label{gamzetdot}
Assume that $m > \frac{7}{2}$, and let $\ww \in C^0([0,T],\L^2(m))$ be
a solution of \reff{SV3}. Then the coefficients $\gamma_i$, $\zeta_{ij}$
defined by \reff{gamzetdef} satisfy, for all $\tau \in [0,T]$,
\begin{eqnarray}\nonumber
  \dot \gamma_i(\tau) &=& -\frac{3}{2}\gamma_i(\tau)\ , \quad i = 1,2,3\ ,\\
  \label{gamzeteq}
  \dot \zeta_{ii}(\tau) &=& -\frac{3}{2}\zeta_{ii}(\tau) + \frac{1}{2}
    \int_{\real^3}(v_3(\xi,\tau)^2 -v_i(\xi,\tau)^2)\d\xi\ , \quad i=1,2\ ,\\
  \nonumber
  \dot \zeta_{ij}(\tau) &=& -\frac{3}{2}\zeta_{ij}(\tau) -\int_{\real^3}
   v_i(\xi,\tau)v_j(\xi,\tau)\d\xi\ , \quad 1\le i < j\le 3\ ,
\end{eqnarray}
where $\vv = (v_1,v_2,v_3)$ is the velocity field obtained from $\ww$
via the Biot-Savart law \reff{SBS3}. 
\end{lemma}

\proof We proceed as in the proof of Lemma~\ref{3betadot}. 
Differentiating \reff{gamzetdef} and integrating by parts, we find
\begin{eqnarray*}
  \dot \gamma_i &=& \int_{\real^3} \qq_i\cdot\bigl(\Lambda \ww
    +\rot(\vv\wedge\ww)\bigr)\d\xi  
  = -\frac{3}{2}\gamma_i + \int_{\real^3} \pp_i \cdot \bigl(\frac{1}{2}
    \nabla|\vv|^2 - (\vv\cdot\nabla)\vv\bigr)\d\xi \\
  &=& -\frac{3}{2}\gamma_i + \int_{\real^3} \vv\cdot\bigl((\vv\cdot\nabla)
    \pp_i\bigr)\d\xi = -\frac{3}{2}\gamma_i\ ,
\end{eqnarray*}
because $(\vv\cdot\nabla)\pp_i = \frac{1}{2}\ee_i\wedge\vv \perp \vv$. 
Similarly, for all $(ij) \in S$, we find
$$
  \dot \zeta_{ij} = -\frac{3}{2}\zeta_{ij} + \int_{\real^3} \vv\cdot
  \bigl((\vv\cdot\nabla)\rot \rr_{ij}\bigr)\d\xi\ .
$$
But $\rot \rr_{ii} = \frac{1}{2}(\xi_3\ee_3-\xi_i\ee_i)$ and 
$\rot \rr_{ij} = -\frac{1}{2}(\xi_i\ee_j +\xi_j\ee_i)$ for $i\neq j$, 
hence
$$
  \vv\cdot\bigl((\vv\cdot\nabla)\rot \rr_{ii}\bigr) = \frac{1}{2}
  (v_3^2-v_i^2)\ , \quad \vv\cdot\bigl((\vv\cdot\nabla)\rot \rr_{ij}\bigr)
  = -v_i v_j \quad (i\neq j)\ .
$$
This concludes the proof. \QED

\medskip
The remainder $\RR$ in \reff{decompp} satisfies the equation
$$
  \partial_\tau \RR = \Lambda \RR + Q_2\bigl((\ww\cdot\nabla)\vv - 
  (\vv\cdot\nabla)\ww\bigr)\ , \quad \div \RR = 0\ ,
$$
where $Q_2 : \L^2(m) \to \cW_2$ is the spectral projection (for the 
operator $\Lambda$) onto the subspace $\cW_2$, see Appendix~\ref{spectrum}. 
Our next result describes the second order asymptotics of $\ww(\tau)$ as 
$\tau \to \infty$. 

\begin{theorem}\label{secondasym}
Let $\frac{3}{2} < \nu < 2$, $m > 2\nu+\frac{1}{2}$, and let 
$\ww \in C^0([0,\infty),\L^2(m))$ be the solution of \reff{SV3}
given by Theorem~\ref{zeroasym} with $\mu = 1$. Then there exist 
constants $b_i$, $c_i$, $d_{ij}$, and $K_5$ such that 
$\|\ww(\tau)-\ww_{app}(\tau)\|_m \le K_5\e^{-\nu\tau}\|\ww_0\|_m$ for 
all $\tau \ge 0$, where
\begin{equation}\label{wappdef}
  \ww_{app}(\xi,\tau) = \sum_{i=1}^3 b_i\e^{-\tau}\ff_i(\xi) + 
  \sum_{i=1}^3 c_i \e^{-\frac{3}{2}\tau}\gg_i(\xi) + \sum_{(ij)\in S} 
  d_{ij}\e^{-\frac{3}{2}\tau}\hh_{ij}(\xi)\ .
\end{equation}
\end{theorem}

\proof If $\ww \in C^0([0,\infty),\L^2(m))$ is the solution of \reff{SV3}
given by Theorem~\ref{zeroasym} and $\vv$ is the corresponding velocity
field, we define $\beta_i$, $\gamma_i$, $\zeta_{ij}$, and $\RR$ by 
\reff{betadef}, \reff{decompp}, \reff{gamzetdef}. It is clear that
$\beta_i(\tau) = b_i \e^{-\tau}$ and $\gamma_i(\tau) = c_i 
\e^{-\frac{3}{2}\tau}$, where
$$
   b_i = \int_{\real^3} \pp_i(\xi)\cdot\ww_0(\xi)\d\xi\ , \quad
   c_i = \int_{\real^3} \qq_i(\xi)\cdot\ww_0(\xi)\d\xi\ , \quad
   i = 1,2,3\ .
$$
By Corollary~\ref{zerolp}, there exists $C_1 > 0$ such that $|\vv(\tau)|_2 
\le C_1 \e^{-\tau}\|\ww_0\|_m$ for all $\tau \ge 0$. Thus, it follows 
easily from Lemma~\ref{gamzetdot} that $|\zeta_{ij}(\tau) - d_{ij}
\e^{-\frac{3}{2}\tau}| \le C_2 \e^{-2\tau}\|\ww_0\|_m^2$ for all 
$\tau \ge 0$, where
\begin{eqnarray}\label{3dijform}
  d_{ii} &=& \zeta_{ii}(0) + \frac{1}{2}\int_0^\infty \e^{\frac{3}{2}\tau}
   \int_{\real^3}(v_3(\xi,\tau)^2 -v_i(\xi,\tau)^2)\d\xi\d\tau\ , \quad
   i = 1,2\ ,\\ \nonumber
  d_{ij} &=& \zeta_{ij}(0) - \int_0^\infty \e^{\frac{3}{2}\tau}
   \int_{\real^3} v_i(\xi,\tau)v_j(\xi,\tau)\d\xi\d\tau\ , \quad
   1 \le i < j \le 3\ .
\end{eqnarray}
To bound the remainder $\RR$, we proceed exactly as in the proof of 
Theorem~\ref{firstasym}. By Proposition~\ref{Sgestim}, there exists
$C_3 > 0$ such that $\|\e^{\tau\Lambda}\RR_0\|_m \le C_3 \e^{-\nu\tau}
\|\ww_0\|_m$ for all $\tau \ge 0$. Using the integral equation
$$
   \RR(\tau) = \e^{\tau\Lambda}\RR_0 + \int_0^\tau Q_2 
   \PPhi(\tau-s,\ww(s),\vv(s))\d s\ , 
$$
together with the bound $\|\ww(\tau)\|_m \le K_0 \e^{-\tau}\|\ww_0\|_m$
given by Theorem~\ref{zeroasym}, it is easy to show that $\|\RR(\tau)\|_m 
\le C_4 \e^{-\nu\tau} \|\ww_0\|_m$ for all $\tau \ge 0$. \QED

\begin{corollary}\label{secondlp} Let $\frac{3}{2} < \nu < 2$ and
$m > 2\nu+\frac{1}{2}$. Let $\ww \in C^0([0,\infty),\L^2(m))$ 
be the solution of \reff{SV3} given by Theorem~\ref{zeroasym} with
$\mu = 1$, and let $\vv$ be the corresponding velocity field. There exists 
$K_6 > 0$ such that, for all $\tau \ge 1$, 
\begin{eqnarray}\label{estttA}
  &&|\ww(\tau) - \ww_{app}(\tau)|_p \le K_6 \e^{-\nu\tau}\|\ww_0\|_m\ , 
  \quad 1 \le p \le \infty\ ,\\ \label{estttB}
  &&|\vv(\tau) - \vv_{app}(\tau)|_q \le K_6 
  \e^{-\nu\tau}\|\ww_0\|_m\ , \quad 1 \le q \le \infty\ ,
\end{eqnarray}
where
$$
  \vv_{app}(\xi,\tau) = \sum_{i=1}^3 b_i\e^{-\tau}\vv^{\ff_i}(\xi) + 
  \sum_{i=1}^3 c_i \e^{-\frac{3}{2}\tau}\vv^{\gg_i}(\xi) + \sum_{(ij)\in S} 
  d_{ij}\e^{-\frac{3}{2}\tau}\vv^{\hh_{ij}}(\xi)\ .
$$
Moreover, the bounds \reff{estttA}, \reff{estttB} hold for all $\tau \ge 0$ 
if $p \le 2$, $q \le 6$.
\end{corollary}

In terms of the original variables, Corollary~\ref{secondlp} shows that, 
for all $t \ge 1$, 
\begin{eqnarray}\nonumber
  &&|\oomega(t) - \oomega_{app}(t)|_p \le C t^{-1-\nu+\frac{3}{2p}}
  \|\oomega_0\|_m\ , \quad 1 \le p \le \infty\ ,\\ \label{uappbd}
  &&|\uu(t) - \uu_{app}(t)|_q \le C t^{-\frac{1}{2}-\nu+\frac{3}{2q}}
  \|\oomega_0\|_m\ , \quad 1 \le q \le \infty\ ,
\end{eqnarray}
where $\oomega_{app}(x,t)$, $\uu_{app}(x,t)$ are the self-similar 
vector fields given by
\begin{eqnarray}\nonumber
   \oomega_{app}(x,t) &=& \frac{1}{1+t}\ww_{app}\Bigl(\frac{x}{\sqrt{1+t}},
   \log(1+t)\Bigr)\ , \\ \label{uappdef}
   \uu_{app}(x,t) &=& \frac{1}{\sqrt{1+t}}\vv_{app}\Bigl(\frac{x}{\sqrt{1+t}},
   \log(1+t)\Bigr)\ .
\end{eqnarray}

\begin{remark} In contrast with the two-dimensional case, the second 
order asymptotic expansions of $\oomega(x,t)$ and $\uu(x,t)$ 
contain only integer powers of $(1+t)^{-1/2}$, and not resonant
terms of the form $(1+t)^{-\alpha}\log(1+t)$. However, following 
\cite{gallay:2001}, one can show that such logarithmic terms do appear 
in the third order asymptotics. This is the reason why the case $\nu = 2$ 
is excluded in Theorem~\ref{secondasym}. In fact, if $m > 9/2$, the 
proof of Theorem~\ref{secondasym} yields the estimate
$$
  \|\ww(\tau)-\ww_{app}(\tau)\|_m \le K_5 (1+\tau)\e^{-2\tau}\|\ww_0\|_m\ ,
  \quad \tau \ge 0\ ,
$$
which appears to be optimal. 
\end{remark}

Two prior papers which discuss the second order asymptotics of solutions 
of the Navier-Stokes equations are those of Carpio, \cite{carpio:1996}, 
and Fujigaki and Miyakawa, \cite{fujigaki:2000}. As we demonstrate below, 
the results of Corollary \ref{secondlp} extend the results of these two 
references. The extension results from the fact that by imposing decay 
conditions on the initial velocity field (see the hypotheses of Theorem~0.6 
in \cite{carpio:1996} and equation (1.2) of \cite{fujigaki:2000}) 
certain terms in the approximating velocity field $\uu_{app}$ are forced 
to be zero. Thus, certain solutions of \reff{NS3} of finite energy 
({\it i.e.} of finite $L^2$ norm) whose asymptotics Corollary \ref{secondlp}
allows us to compute are excluded from consideration by the decay conditions 
of \cite{carpio:1996} and \cite{fujigaki:2000}. This is a further reason 
that we feel it is more natural to impose decay conditions on the vorticity 
rather than the velocity. Note that in deriving the higher-order asymptotics 
in \cite{fujigaki:2000} (Theorem 2.2 (ii)) increasingly stringent decay
conditions are imposed on the velocity which results in more and more terms 
in the asymptotics being zero. To compare the results of 
Corollary~\ref{secondlp} with those of the previous references, first note 
that the requirement that $(1+|x|)\vv_0 \in L^1(\real^3)^3$ implies that
that $b_i = 0$ for $i=1,2,3$ by Corollary~\ref{3velint}. Furthermore, 
by Corollary \ref{3velintbis}, the remaining coefficients in $\vv_{app}$ 
satisfy
$$
  c_1=b_{23}\ , \ c_2 = b_{31}\ ,\ c_3 = b_{12}\ ,
$$
and
$$
  d_{11} = \half (c_{33}-c_{11})\ ,\ d_{22} = \half (c_{33} - c_{22})\ ,
  \ d_{12} = - c_{12} \ ,\ d_{23}=-c_{23}\ ,\ d_{13}=-c_{13}\ ,
$$
where $b_{k\ell}$ and $c_{k\ell}$ are defined in \reff{bcmom}.
The expressions above for $d_{ij}$ follow easily from \reff{u-v3}, 
\reff{3dijform} if one remembers that $\zeta_{ij}(0) = 0$ by 
Corollary~\ref{3velintbis}. 

Now consider the term $\AA(\xi) = \sum_{j=1}^3 c_j \vv^{\gg_j}(\xi)$, in
$\vv_{app}$.  Since $\vv^{\gg_j} = \ff_j = \half (\ee_j \wedge \xi)G$, 
where $G(\xi) = (4\pi)^{-3/2} \exp(-|\xi|^2/4)$, this sum can be rewritten
as
$$
  \AA(\xi) = \half \sum_{j=1}^3 c_j (\ee_j \wedge \xi) G(\xi) = \half 
  G(\xi) (\cc \wedge \xi)\ ,
$$
where $\cc = (c_1,c_2,c_3)$. Examining this expression component-by-component,
we see that
$$
  A_1(\xi) = \half G(\xi) (c_2\xi_3 - c_3 \xi_2) = \half G(\xi) 
  (b_{31} \xi_3 - b_{12} \xi_2) = - (b_{31} \partial_3 G + b_{21} 
  \partial_2 G)\ ,
$$
using the anti-symmetry of the $b_{jk}$. The other components are treated 
in like fashion and we find
$$
   A_j(\xi) =  -  \sum_{k=1}^3 b_{k j} \partial_k G(\xi)\ , \quad 
   j = 1,2,3\ .
$$
To treat the term $\BB(\xi) = \sum_{(ij) \in S} d_{ij} \vv^{\hh_{ij}}(\xi)$,
we note that using the expressions for $\vv^{\hh_{ij}}$ in 
Appendix~\ref{velocity3d}, one has
$$
  \vv^{\hh_{ij}} = 2 \partial_i \partial_j (\nabla \Phi) +
  (\partial_i G) \ee_j + (\partial_j G) \ee_i\ ,
$$
where $-\Delta \Phi = G$.  Then using the expressions for $d_{ij}$ from
above a straightforward computation shows that
$$
  \BB(\xi)= - \sum_{i,j=1}^3 c_{ij} \left(\partial_i \partial_j 
  (\nabla \Phi) + (\partial_i G)\ee_j \right)\ .
$$
In \cite{fujigaki:2000}, the following notation is introduced:
$$
   E_t(x) = \frac{1}{t^{3/2}}\,G\Bigl(\frac{x}{\sqrt{t}}\Bigr)\ ,
   \quad F_{\ell,jk}(x,t) = \frac{1}{t^2}\,F_{\ell,jk}
   \Bigl(\frac{x}{\sqrt{t}}\,,\,1\Bigr)\ , \quad x \in \real^3\ ,
   \quad t > 0\ ,
$$
where $F_{\ell,jk}(\cdot,1) = \partial_k \partial_{\ell} (\partial_j \Phi)
+ (\partial_{\ell} G) \delta_{jk}$. Using \reff{uappdef} and the 
expressions above for $\AA(\xi)$ and $\BB(\xi)$, it is easy to 
verify that the velocity field $\uu_{app}$ in \reff{uappbd} satisfies
$$
   \ee_j \cdot \uu_{app}(x,t{-}1) = -\sum_{k=1}^3 b_{kj} \partial_k E_t(x)
   -\sum_{k,\ell=1}^3 c_{k\ell}F_{\ell,jk}(x,t)\ , \quad 
   j = 1,2,3\ .
$$
Thus, written in terms of this notation, Corollary \ref{secondlp}
implies:

\begin{corollary} \label{fmcomp} If in addition to the hypotheses
of Corollary~\ref{secondlp}, one assumes that the initial condition for
the velocity field satisfies $\rho \uu_0 \in L^1(\real^3)^3$, then 
for $j = 1,2,3$ one has the estimate
$$
  \Big|u_j (\cdot,t) + \sum_{k=1}^3 b_{kj} \partial_k E_t(\cdot)
  + \sum_{k,\ell =1}^3 c_{k\ell}F_{\ell,jk}(\cdot,t)\Big|_q
  \le C t^{-\half - \nu +\frac{3}{2q}} \|\oomega_0\|_m\ ,
  \quad t \ge 1\ .
$$
\end{corollary}
Comparing with equation (2.4) of \cite{fujigaki:2000} we see that this is 
compatible with the results of Fujigaki and Miyakawa. Rewriting the 
asymptotics in a slightly different way one finds that they also agree with 
Theorem~0.6 of \cite{carpio:1996}.

%%%%%%%%%%%%%%%%%%%%%%%%%%%%%%%%%%%%%%%%%%%%%%%%%%%%%%%%%%%%%%%%%%%%%%%%%%%

\section{The strong-stable manifold of the origin}\label{ssmanif}

In this section, we assume that $m > 7/2$ and we consider in more detail 
the dynamics of \reff{SV3} in the invariant subspace $\cW_1$ of $\L^2(m)$ 
defined by \reff{Wndef}. If $\ww_0 \in \cW_1$ satisfies $\|\ww_0\|_m \le
r_0$, where $r_0 > 0$ is as in Theorem~\ref{zeroasym}, the solution
$\ww(\cdot,\tau)$ of \reff{SV3} with initial data $\ww_0$ can be 
decomposed as
\begin{equation}\label{3wwvvdecomp}
  \ww(\xi,\tau) = \sum_{i=1}^3 \gamma_i(\tau)\gg_i(\xi) + \sum_{(ij)\in S} 
   \zeta_{ij}(\tau)\hh_{ij}(\xi) + \RR(\xi,\tau)\ ,
\end{equation}
where $\gamma_i$, $\zeta_{ij}$ are defined in \reff{gamzetdef} and 
$\RR(\cdot,\tau)$ belongs to the subspace $\cW_2$ of $\L^2(m)$. As for 
the velocity field, we have
$$
  \vv(\xi,\tau) = \sum_{i=1}^3 \gamma_i(\tau)\vv^{\gg_i}(\xi) + 
  \sum_{(ij)\in S} \zeta_{ij}(\tau)\vv^{\hh_{ij}}(\xi) +\vv^{\RR}(\xi,\tau)\ .
$$
Setting $b_i = 0$ in Theorem~\ref{secondasym}, we see that 
$\|\ww(\cdot,\tau)\|_m = \cO(\e^{-\frac{3}{2}\tau})$ as $\tau \to +\infty$. 
We now define the {\it local strong-stable manifold} of the origin by
\begin{equation}\label{3ssmanif}
  W_s^{loc} = \Bigl\{\ww_0 \in \L^2(m) \,\Big|\, \|\ww_0\|_m \le r_0\,,~
  \lim_{\tau\to\infty} \e^{\frac{3}{2}\tau}\|\Phi_\tau \ww_0\|_m = 0
  \Bigr\}\ ,
\end{equation}
where $\Phi_\tau \ww_0 = \ww(\tau)$ is the solution of \reff{SV3} in 
$\L^2(m)$ with initial data $\ww_0$. It is clear from 
Theorem~\ref{firstasym} that $W_s^{loc} \subset \cW_1$. On the other
hand, it follows from Lemma~\ref{gamzetdot} that $W_s^{loc} \not\subset 
\cW_2$. Indeed, the inclusion $W_s^{loc} \subset \cW_2$ would mean 
that the integrals in the right-hand side of \reff{gamzeteq} vanish 
identically for vorticities $\ww \in \cW_2$, but it is easy to verify
that this is not the case. 

By construction, $\cW_1 = \cW_2 \oplus \cV$, where $\cV$ is the 
eight-dimensional space spanned by the vector fields $\gg_i$ for $i = 1,2,3$ 
and $\hh_{ij}$ for $(ij) \in S$. Using invariant manifold theory as in 
\cite{gallay:2001}, it is rather straightforward to show that $W_s^{loc}$ 
is a smooth submanifold of $\cW_1$ which is tangent at the origin to the 
subspace $\cW_2$. In other words, there exists a smooth function $f : \cW_2 
\to \cV$ satisfying $f(0) = 0$, $f'(0) = 0$, and such that 
$W_s^{loc} = G(f) \cap B(r_0)$, where
$$
   G(f) = \{\ww + f(\ww)\,|\, \ww \in \cW_2\} \subset \cW_1\ ,
   \quad B(r_0) = \{\ww \in \cW_1 \,|\, \|\ww\|_m \le r_0\}\ .
$$
In particular, the manifold $W_s^{loc}$ is of codimension $8$ in 
$\cW_1$, hence of codimension $11$ in $\L^2(m)$. By construction, 
$W_s^{loc}$ is locally positively invariant in the following sense:
if $\ww_0 \in W_s^{loc}$, then $\ww(\tau) = \Phi_\tau \ww_0 \in W_s^{loc}$
as long as $\ww(\tau) \in B(r_0)$. If in addition $\|\ww_0\|_m \le r_0/K_0$, 
where $K_0$ is as in Theorem~\ref{zeroasym}, then $\ww(\tau) \in W_s^{loc}$
for all $\tau \ge 0$ and $\e^{\frac{3}{2}\tau}\|\ww(\tau)\|_m \to 0$ as
$\tau \to +\infty$. 

\figurewithtex 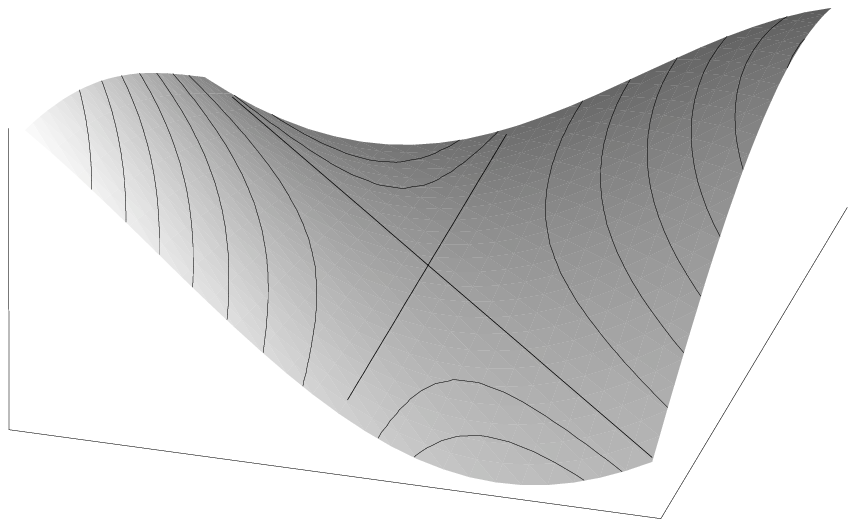 Fig2.tex 6.0 9.0 {\bf Fig.~2:} A {\it schematic}
picture of the local strong-stable manifold $W_s^{loc}$ (shaded
surface). The horizontal plane is the infinite-dimensional
subspace $\cW_2$, and the vertical axis the eight-dimensional 
space $\cV$. The intersection $W_s^{loc} \cap \cW_2$, which is 
also infinite-dimensional, is represented by two line segments.\cr

\medskip
The following result is characterization of the local strong-stable manifold. 

\begin{proposition}\label{3caractws}
Fix $m > 7/2$, and assume that $\ww_0 \in \cW_1 \subset \L^2(m)$ satisfies 
$\|\ww_0\|_m \le r_0$, where $r_0 > 0$ is as in Theorem~\ref{zeroasym}. 
Let $w(\xi,\tau)$ be the solution of \reff{SV3} with initial data $\ww_0$, 
and let $\vv(\xi,\tau)$ be the corresponding velocity field. Define
the functions $\gamma_i(\tau), \zeta_{ij}(\tau)$ by \reff{gamzetdef} 
and the coefficients $c_{k\ell}$ by \reff{bcmom}, with
\begin{equation}\label{3uvvscaling}
  \uu(x,t) = \frac{1}{\sqrt{1+t}}\,\vv\Bigl(\frac{x}{\sqrt{1+t}},
  \log(1+t)\Bigr)\ .
\end{equation}
Then the following three statements are equivalent:\\
\null\quad {\rm 1)} $\displaystyle{\lim_{\tau\to\infty} \e^{\frac{3}{2}\tau} 
  \|\ww(\cdot,\tau)\|_m = 0}$, namely $\ww_0 \in W_s^{loc}$.\\
\null\quad {\rm 2)} $\displaystyle{\lim_{t\to\infty} t^{5/4}|\uu(\cdot,t)|_2 
  = 0}$.\\
\null\quad {\rm 3)} $\gamma_1(0) = \gamma_2(0) = \gamma_3(0) = 0$,\\
\null\qquad $\zeta_{11}(0) = \frac{1}{2}(c_{11}-c_{33})$, 
  $\zeta_{22}(0) = \frac{1}{2}(c_{22}-c_{33})$, $\zeta_{12}(0) = c_{12}$, 
  $\zeta_{13}(0) = c_{13}$, $\zeta_{23}(0) = c_{23}$.
\end{proposition}

\proof We apply Theorem~\ref{secondasym} with $\frac{3}{2} < \nu < 
\frac{m}{2}-\frac{1}{4}$. Since $\ww_0 \in \cW_1$, we have $b_i = 0$
for $i = 1,2,3$, so that $\|\ww(\tau)\|_m = \cO(\e^{-\frac{3}{2}\tau})$ as 
$\tau \to +\infty$. We shall show that statements 1), 2), 3) in 
Proposition~\ref{3caractws} are all equivalent to\\
\null\quad 4) $c_i = 0$ for $i = 1,2,3$ and $d_{ij} = 0$ for all $(ij) 
\in S$.\\
Indeed, it is clear from \reff{wappdef} that $\|\ww_{app}(\tau)\|_m = 
K\e^{-\frac{3}{2}\tau}$, where $K = 0$ if and only if 4) holds. Thus 
${\rm 1)} \Leftrightarrow {\rm 4)}$ by Theorem~\ref{secondasym}. Similarly, 
it follows from \reff{uappdef} that $|\uu_{app}(t)|_2 = K' (1+t)^{-5/4}$, 
where $K' = 0$ if and only if 4) holds. Since $|\uu(t)-\uu_{app}(t)|_2 \le
C(1+t)^{\frac{1}{4}-\nu}$ by \reff{uappbd}, we conclude that ${\rm 2)} 
\Leftrightarrow {\rm 4)}$. Finally, using \reff{bcmom}, \reff{3dijform}
and the change of variables \reff{3uvvscaling}, we obtain the 
relations
\begin{eqnarray*}
   d_{ii} &=& \zeta_{ii}(0) + \frac{1}{2}(c_{33}-c_{ii})\ , \quad 
     i = 1,2\ ,\\
   d_{ij} &=& \zeta_{ij}(0) -c_{ij}\ , \quad 1 \le i < j \le 3\ .
\end{eqnarray*}
We also know that $c_i = \gamma_i(0)$ for $i = 1,2,3$. Therefore, 
${\rm 3)} \Leftrightarrow {\rm 4)}$. \QED

\medskip
Using Proposition~\ref{3caractws}, it is easy to prove Theorem~\ref{MiSchon}
in the case where $\ww_0 = \rot(\uu_0)$ satisfies $\|\ww_0\|_m \le r_0$ for 
some $m > 7/2$. Indeed, let $\ww(\xi,\tau)$ be the solution of \reff{SV3} 
with initial data $\ww_0$, and define $\gamma_i(\tau), \zeta_{ij}(\tau)$
by \reff{gamzetdef}. By Corollary~\ref{3velintbis}, the assumption 
$\rho \uu_0 \in L^1(\real^3)^3$ implies that $\ww_0 \in \cW_1$ (so that
$\ww(\tau) \in \cW_1$ for all $\tau \ge 0$) and that $\zeta_{ij}(0) = 0$ 
for all $(ij) \in S$. 

\smallskip\noindent{\bf a)}
Assume first that \reff{optidecay} holds, namely $\ww_0 \in W_s^{loc}$. 
Then point 3) in Proposition~\ref{3caractws} shows that $c_{11} = 
c_{22} = c_{33}$ and $c_{ij} = 0$ if $i \neq j$, hence the matrix 
$(c_{k\ell})$ is scalar. In addition, since $\gamma_i(0) = 0$, we have 
$b_{k\ell} = 0$ by Corollary~\ref{3velintbis}. This proves 
\reff{more-moments}. 

\smallskip\noindent{\bf b)}
Conversely, assume that \reff{more-moments} holds. Then $\gamma_i(0) = 0$
by Corollary~\ref{3velintbis}, and since $\zeta_{ii}(0) = 0 = 
\frac{1}{2}(c_{ii}-c_{33})$ for $i = 1,2$, $\zeta_{ij}(0) = 0 = c_{ij}$ 
for $i \neq j$, it follows from Proposition~\ref{3caractws} that 
$\ww_0 \in W_s^{loc}$. This concludes the argument. \QED

\medskip
As is clear from this proof, if $\ww_0$ lies in $W_s^{loc}$ and 
if the corresponding velocity field $\vv_0$ satisfies $\rho \vv_0 \in 
L^1(\real^3)^3$, then necessarily $\ww_0 \in \cW_2$. Thus, from our
point of view, Theorem~\ref{MiSchon} is a characterization of the 
noninvariant set $W_s^{loc} \cap \cW_2$, and not of $W_s^{loc}$ itself. 
As was already observed, $W_s^{loc}$ is not contained in $\cW_2$, 
and it is not a priori obvious that $W_s^{loc} \cap \cW_2 \neq \{0\}$! 
In fact, using a nice argument due to L. Brandolese, it turns out that
$W_s^{loc} \cap \cW_2$ is infinite-dimensional. 
Following \cite{brandolese:2001}, we say that a vector field 
$\uu : \real^3 \to \real^3$ is {\it symmetric} if it satisfies
the following two properties:\\
A) $u_1(x_1,x_2,x_3) = u_2(x_3,x_1,x_2) = u_3(x_2,x_3,x_1)$ for all 
   $x = (x_1,x_2,x_3) \in \real^3$.\\
B) For all $i \in \{1,2,3\}$, $u_i(x_1,x_2,x_3)$ is an {\it odd} function
   of $x_i$ and an {\it even} function of $x_j$\\
\null\quad\ for all $j \neq i$.\\
If $\uu$ is symmetric, then $\Delta \uu$ and $(\uu\cdot\nabla)\uu$ are
also symmetric. This implies, roughly speaking, that the space of 
symmetric velocity fields is invariant under the Navier-Stokes evolution
(whenever defined). 

Assume now that $\vv : \real^3 \to \real^3$ is symmetric and that
the vorticity $\ww = \rot \vv$ belongs to $\L^2(m)$ for some $m > 7/2$. 
Then $\ww$ satisfies:\\
A')  $w_1(\xi_1,\xi_2,\xi_3) = w_2(\xi_3,\xi_1,\xi_2) = 
    w_3(\xi_2,\xi_3,\xi_1)$ for all $\xi \in \real^3$.\\
B') For all $i \in \{1,2,3\}$, $w_i(\xi_1,\xi_2,\xi_3)$ is an {\it even} 
    function of $\xi_i$ and an {\it odd} function of $\xi_j$\\
\null\quad\ for all $j \neq i$.\\
Using these properties together with \reff{betadef}, \reff{gamzetdef},
we discover that $\beta_i = \gamma_i = \zeta_{ij} = 0$ for all $i,j$, 
hence $\ww \in \cW_2$. On the other hand, it is clear that the integrals
in the right-hand side of \reff{gamzeteq} vanish identically if 
$\vv$ is symmetric, so that $c_{k\ell} = 0$ in \reff{bcmom}. Thus, if 
$\|\ww\|_m \le r_0$, it follows from Proposition~\ref{3caractws} that 
$\ww \in W_s^{loc}$. Summarizing, we have shown:
$$
   W_s^{loc} \cap \cW_2 \supset \Bigl\{\ww \in \L^2(m)\,\Big|\, 
   \|\ww\|_m \le r_0\ ,~ \ww = \rot \vv \hbox{ with } \vv 
   \hbox{ symmetric}\Bigr\}\ .
$$

We conclude this section with a somewhat surprising observation. Let 
$\Psi_t$ be the local semiflow defined by the vorticity equation \reff{V3} 
in $\L^2(m)$ for $m > 7/2$. If $\oomega \in \L^2(m)$, let $\UU_\omega$
be the velocity field obtained from $\oomega$ via the Biot-Savart law
\reff{BS3}. Then
$$
   W_s^{loc} = \Bigl\{\oomega_0 \in \L^2(m)\,\Big|\, \|\oomega_0\|_m 
   \le r_0\,,~ \lim_{t\to+\infty}t^{5/4}|\UU_{\Psi_t \omega_0}|_2 = 0
   \Bigr\}\ .
$$
This characterization follows from the equivalence ${\rm 1)} 
\Leftrightarrow {\rm 2)}$ in Proposition~\ref{3caractws} and from the 
fact that the change of variables \reff{omega-w3}, \reff{u-v3} reduces 
to the identity when $t = 0$. As a consequence, $W_s^{loc}$ is locally 
invariant under {\it both} semiflows $\Psi_t$ and $\Phi_\tau$, although 
the orbits of the same initial point under $\Psi_t$ and $\Phi_\tau$ are 
of course different! This curious property originates in the fact that 
in both the original and rescaled variables the strong-stable manifold 
can be characterized in terms of the decay rate (as time goes to infinity)
of the solutions lying in it, see \cite{gallay:2001} for a more detailed 
discussion. In concrete terms, the observation above implies that 
the picture of $W_s^{loc}$ in Fig.~2 is not affected at all when we
return to the original variables using \reff{omega-w3}. 

%%%%%%%%%%%%%%%%%%%%%%%%%%%%%%%%%%%%%%%%%%%%%%%%%%%%%%%%%%%%%%%%%%%%%%%%%%%
\appendix
\section{Spectrum of the operator $\Lambda$}\label{spectrum}

In Appendix~A of \cite{gallay:2001}, we study in detail the linear
operator
\begin{equation}\label{cLop}
  \cL = \Delta_\xi + \frac{1}{2}\xi \cdot \nabla_\xi + \frac{N}{2}\ ,
  \quad \xi \in \real^N\ , \quad N \ge 1\ ,
\end{equation}
acting on the function $L^2(m) = \{f \in L^2(\real^N) \,|\, \|f\|_m < 
\infty\}$, where
$$
   \|f\|_m = \left(\int_{\real^N}(1+|\xi|)^{2m} |f(\xi)|^2 \d\xi 
   \right)^{1/2}\ .
$$
In particular, we determine exactly the spectrum of $\cL$:

\begin{theorem}\label{Lspectrum} {\rm \cite{gallay:2001}}
Fix $m \ge 0$, and let $\cL$ be the linear operator \reff{cLop} in 
$L^2(m)$, defined on its maximal domain. Then the spectrum of 
$\cL$ is
$$
\sigma(\cL) = \Bigl\{\lambda \in \complex\,\Big|\, \Re(\lambda) \le 
\frac{N}{4} - \frac{m}{2}\Bigr\} \cup \Bigl\{-\frac{k}{2} \,\Big|\, 
k \in \intplus\Bigr\}\ .
$$
Moreover, if $m > \frac{N}{2}$ and if $k \in \intplus$ satisfies 
$k + \frac{N}{2} < m$, then $\mu_k = -\frac{k}{2}$ is an isolated 
eigenvalue of $\cL$, with multiplicity ${N+k-1 \choose k}$.
\end{theorem}

The eigenfunctions corresponding to the isolated eigenvalues 
$\mu_k = -\frac{k}{2}$ can be computed explicitly. Moreover, it is 
shown in \cite{gallay:2001} that $\cL$ generates a $C_0$ semigroup 
$\e^{\tau\cL}$ in $L^2(m)$, and sharp estimates are obtained for the 
norm of $\e^{\tau\cL}$ in various spectral subspaces of $L^2(m)$.

\medskip
In this section, we adapt the results in \cite{gallay:2001} to the 
particular case where $N = 3$ and $\cL \equiv \Lambda +\frac{1}{2}$ acts
on the space of {\it divergence free} vector fields $\L^2(m)$ defined
in \reff{LL2m}. Remark that $\div(\Lambda \ff) = \cL \div(\ff)$, so 
that $\Lambda$ preserves the divergence free condition. The analogue of 
Theorem~\ref{Lspectrum} is:

\begin{theorem}\label{Lamspect}
Fix $m \ge 0$, and let $\Lambda$ be the linear operator \reff{Lamop} in 
$\L^2(m)$, defined on its maximal domain. Then the spectrum of 
$\Lambda$ is
$$
  \sigma(\Lambda) = \Bigl\{\lambda \in \complex\,\Big|\, \Re(\lambda) \le 
  \frac{1}{4} - \frac{m}{2}\Bigr\} \cup \Bigl\{-\frac{k+1}{2} \,\Big|\, 
  k \in \intplus^*\Bigr\}\ ,
$$
where $\intplus^* = \intplus\setminus\{0\}$. Moreover, if $m > \frac{5}{2}$ 
and if $k \in \intplus^*$ satisfies $k + \frac{3}{2} < m$, then 
$\lambda_k = -\frac{k+1}{2}$ is an isolated eigenvalue of $\cL$, with 
multiplicity $k(k+2)$.
\end{theorem}

\proof We first discuss the {\sl discrete} spectrum of $\Lambda$. 
Fix $k \in \intplus$, and take $\alpha = (\alpha_1,\alpha_2,\alpha_3) 
\in \intplus^3$ such that $|\alpha| = \alpha_1 + \alpha_2 + \alpha_3 = k$.
Then the Hermite function $\phi_\alpha : \real^3 \to \real$ defined by
$$
   \phi_\alpha = \partial^\alpha G \equiv \partial_1^{\alpha_1}
   \partial_2^{\alpha_2} \partial_3^{\alpha_3} G\ , \quad
   \hbox{ where }G(\xi) = \frac{1}{(4\pi)^{3/2}} \,\e^{-|\xi|^2/4}\ ,
$$
is an eigenfunction of $\cL$ with eigenvalue $-\frac{k}{2}$. Let
$E_k = \span\{\phi_\alpha \,|\, \alpha \in \intplus^3\ ,~|\alpha| = k\}$
and
$$
  \cE_k = \{\ff = (f_1,f_2,f_3) \,|\, \div \ff = 0\ ,~
   f_i \in E_k \hbox{ for } i = 1,2,3\}\ .
$$
By construction, $\cE_k \subset \L^2(m)$ for all $m \ge 0$, and any $\ff \in
\cE_k$ satisfies $\Lambda \ff = -\frac{k+1}{2}\ff$. Moreover, using
the characterization of $\cE_k$ in Fourier variables (see \cite{gallay:2001}), 
it is not difficult to show that $\dim(\cE_k) = k(k+2)$. In particular, 
for any $k \in \intplus^*$, $\lambda_k = -\frac{k+1}{2}$ is an 
eigenvalue of $\Lambda$ with multiplicity (at least) $k(k+2)$. 

\smallskip
Next, fix $\lambda \in \complex$ such that $\Re(\lambda) < 1/4$ 
and $-(\lambda+1)\notin \intplus$. Proceeding as in \cite{gallay:2001}, 
it is easy to verify that the function $\psi_\lambda : \real^3 \to \real$
defined in Fourier variables by
$$
  \hat \psi_\lambda(p) = |p|^{-2(\lambda+1)} \e^{-|p|^2}
  (-\I p_2, \I p_1, 0)
$$
satisfies $\Lambda \psi_\lambda = \lambda \psi_\lambda$ and 
$\div \psi_\lambda = 0$. Moreover, $\psi_\lambda \in \L^2(m)$ if and
only if $\Re(\lambda) < \frac{1}{4}-\frac{m}{2}$. This shows that
$\sigma(\Lambda) \supset \{\lambda \in \complex\,|\, \Re(\lambda) \le 
\frac{1}{4} - \frac{m}{2}\}$. 

\smallskip
Now, fix $n \in \allint$ and assume that $m \ge 0$, $m > n+\frac{3}{2}$. 
Let
$$
   W_n = \Bigl\{f \in L^2(m) \,\Big|\, \int_{\real^3} 
   \xi^\alpha f(\xi)\d\xi = 0 \hbox{ for all }\alpha \in \intplus^3
   \hbox{ with }|\alpha| \le n\Bigr\}\ .
$$
In particular, $W_n = L^2(m)$ if $n < 0$. We define closed subspaces 
$\cV_n$, $\cW_n$ of $\L^2(m)$ by $\cV_n = \oplus_{k=1}^n\cE_k$ and
\begin{equation}\label{Wndef}
   \cW_n = \bigl\{\ff \in \L^2(m) \,\big|\, f_i \in W_n 
   \hbox{ for }i=1,2,3 \bigr\}\ .
\end{equation}
By definition, $\cV_n = \{0\}$ and $\cW_n = \L^2(m)$ if $n \le 0$. 
(We recall that any $\ff \in \L^2(m)$ with $m > \frac{3}{2}$ satisfies 
$\int_{\real^3}\ff(\xi)\d\xi = 0$ as a consequence of the divergence free 
condition; hence $\cW_0 = \L^2(m)$.) Using again the characterization of 
$\L^2(m)$ in Fourier variables, it is easy to verify that $\L^2(m) = 
\cV_n \oplus \cW_n$. Let $P_n : \L^2(m) \to \L^2(m)$ be the (unique) 
continuous projection satisfying $\range(P_n) = \cV_n$, $\ker(P_n) = 
\cW_n$, and let $Q_n = \oone -P_n$. In particular, $P_n = 0$ and 
$Q_n = \oone$ for all $n \in \allint$, $n \le 0$. The following estimates 
on the semigroup $\e^{\tau\Lambda} = \e^{-\tau/2}\e^{\tau\cL}$ are proved
in \cite{gallay:2001} (Proposition~A.2):

\begin{proposition}\label{Sgestim} Let $a(\tau) = 1-\e^{-\tau}$, 
$\tau \ge 0$.\\
{\bf (a)} Fix $m \ge 0$, and take $n \in \allint$ such that
$n + \frac{3}{2} < m \le n+\frac{5}{2}$. For all $\alpha \in \intplus^3$
and all $\epsilon > 0$, there exists $C > 0$ such that
\begin{equation}\label{Sg1}
  \|\partial^\alpha \e^{\tau\Lambda}f\|_m \le \frac{C}{a(\tau)^{|\alpha|/2}}
  \,\e^{\frac{\tau}{2}(\frac{1}{2}-m+\epsilon)} \|f\|_m\ , \quad
  \tau > 0\ ,
\end{equation}
for all $f \in W_n \subset L^2(m)$. 

\smallskip\noindent
{\bf (b)} Fix $n \in \intplus \cup \{-1\}$, and take $m \in \real$ such that
$m > n+\frac{5}{2}$. For all $\alpha \in \intplus^3$ and all $\epsilon > 0$, 
there exists $C > 0$ such that
\begin{equation}\label{Sg2}
  \|\partial^\alpha \e^{\tau\Lambda}f\|_m \le \frac{C}{a(\tau)^{|\alpha|/2}}
  \,\e^{-\frac{n+2}{2}\tau} \|f\|_m\ , \quad \tau > 0\ ,
\end{equation}
for all $f \in W_n \subset L^2(m)$.
\end{proposition}

If $m$ and $n$ are as in part {\bf (a)} of Proposition~\ref{Sgestim}, 
it follows from \reff{Sg1} that
$$
   \|\e^{\tau\Lambda}Q_n \ff\|_m \le C \e^{\frac{\tau}{2}(\frac{1}{2}-m
   +\epsilon)} \|\ff\|_m\ , \quad \tau \ge 0\ ,
$$
for all $\ff \in \L^2(m)$. By the Hille-Yosida theorem, this implies
that $\sigma(\Lambda Q_n) \subset \{\lambda \in \complex\,|\, \Re(\lambda) 
\le \frac{1}{4}-\frac{m}{2}\}$. On the other hand, by construction, 
we have $\sigma(\Lambda P_n) = \emptyset$ if $n \le 0$ and 
$\sigma(\Lambda P_n) = \{-1;-\frac{3}{2};\dots;-\frac{n+1}{2}\}$ if 
$n \in \intplus^*$. In particular, $\sigma(\Lambda P_n) \cap 
\sigma(\Lambda Q_n) = \emptyset$, hence the multiplicity of the eigenvalue 
$\lambda_k$ ($k = 1,\dots,n$) is exactly $k(k+2)$. Finally, since
$\sigma(\Lambda) \subset \sigma(\Lambda P_n) \cup \sigma(\Lambda Q_n)$, 
we have
$$
  \sigma(\Lambda) \subset \Bigl\{\lambda \in \complex\,\Big|\, \Re(\lambda) 
  \le \frac{1}{4} - \frac{m}{2}\Bigr\} \cup \Bigl\{-\frac{k+1}{2} \,\Big|\, 
  k \in \intplus^*\Bigr\}\ .
$$
This concludes the proof of Theorem~\ref{Lamspect}.\QED 

\medskip
The estimates in Proposition~\ref{Sgestim} can be generalized 
to weighted $L^p$ spaces with $p \neq 2$. For our purposes in this paper, 
the following result will be sufficient:

\begin{proposition}\label{LpLqestim} Let $1 \le q \le p \le \infty$,
$m \ge 0$ and $T > 0$. For all $\alpha \in \intplus^3$, there
exists $C > 0$ such that, for all $f \in L^q(m)$,
\begin{equation}\label{LpLq}
  |\rho^m \partial^\alpha \e^{\tau\Lambda}f|_p \le 
  \frac{C}{a(\tau)^{\frac{3}{2} (\frac{1}{q}-\frac{1}{p})
  +\frac{|\alpha|}{2}}} \,|\rho^m f|_q\ ,
  \quad 0 < \tau \le T\ ,
\end{equation}
where $\rho(\xi) = 1+|\xi|$. 
\end{proposition}

\proof See \cite{gallay:2001}, Proposition~A.5. 

%%%%%%%%%%%%%%%%%%%%%%%%%%%%%%%%%%%%%%%%%%%%%%%%%%%%%%%%%%%%%%%%%%%%%%%%%%%

\section{Bounds on the velocity field}
\label{velocity3d}

We first list a few identities that are satisfied by the vorticity 
$\ww$ as a consequence of the divergence free condition. If $\ww \in
\L^1(\real^3)^3$, then
\begin{equation}\label{Aid}
  \int_{\real^3} w_i(\xi)\d\xi = 0 \quad \hbox{for all } i \in \{1,2,3\}\ ,
\end{equation}
because $w_i = \div(\xi_i \ww)$. This fact is not hard to prove in Fourier 
variables, but has been overlooked in some papers on the subject until 
recently. If in addition $\rho \ww \in L^1(\real^3)^3$ (where
$\rho(\xi) = 1+|\xi|$), then
\begin{equation}\label{Bid}
  \int_{\real^3} (\xi_i w_j(\xi) + \xi_j w_i(\xi))\d\xi = 0 \quad 
  \hbox{for all } i,j \in \{1,2,3\}\ ,
\end{equation}
because $\xi_i w_j + \xi_j w_i = \div(\xi_i \xi_j \ww)$. Thus only
three first-order moments of $\ww$ (out of nine) are possibly nonzero. 
Finally, if we assume that $\rho^2 \ww \in L^1(\real^3)^3$, then
\begin{equation}\label{Cid}
  \int_{\real^3} (\xi_i \xi_j w_k(\xi) + \xi_j \xi_k w_i(\xi) + 
  \xi_k \xi_i w_j(\xi)) \d\xi = 0 \quad 
  \hbox{for all } i,j,k \in \{1,2,3\}\ ,
\end{equation}
because $\xi_i \xi_j w_k + \xi_j \xi_k w_i + \xi_k \xi_i w_j = 
\div(\xi_i \xi_j \xi_k \ww)$. This means that only eight second-order
moments of $\ww$ (out of eighteen) are possibly nonzero. 

\medskip
Next, we give explicit formulas for the velocity fields corresponding, 
via the Biot-Savart law \reff{BS3}, to the first eigenfunctions 
of the linear operator $\Lambda$ acting on $\L^2(m)$.  

\smallskip\noindent
{\bf 1)} If $m > 5/2$, the first eigenvalue $\lambda_2 = -1$ has 
multiplicity three. A basis of eigenfunctions is $\{\ff_1,\ff_2,\ff_3\}$, 
where $\ff_i = \pp_i G$ and $\pp_1,\pp_2,\pp_3$ are defined in \reff{ppdef}. 
The corresponding velocity fields are
\begin{equation}\label{vfidef}
  \vv^{\ff_i} = \partial_i (\nabla\Phi) + G\ee_i\ , \quad i = 1,2,3\ ,
\end{equation}
where
$$
  \Phi(\xi) = \frac{1}{(4\pi)^{3/2}} \frac{2}{|\xi|}
  \int_0^{|\xi|} \e^{-z^2/4}\d z\ , \quad \xi \in \real^3\setminus\{0\}\ .
$$
Remark that $-\Delta\Phi = G$, so that $\div \vv^{\ff_i} = 0$ and $\rot 
\vv^{\ff_i} = \rot(G\ee_i) = \ff_i$. A direct calculation shows that
$\vv^{\ff_i} \sim |\xi|^{-3}$ as $|\xi| \to \infty$, hence $\vv^{\ff_i}
\notin \L^1(\real^3)$. 

\smallskip\noindent
{\bf 2)} If $m > 7/2$, the second eigenvalue $\lambda_3 = -3/2$ has 
multiplicity eight. A convenient basis of eigenfunctions is given by 
$\{\gg_i\}_{i=1,2,3}$ and $\{\hh_{ij}\}_{(ij)\in S}$, where $\gg_i$ 
are defined in \reff{ggdef} and $\hh_{ij}$ in \reff{hhdef}. (We recall that 
$S = \{(11),(12),(13),(22),(23)\}$.) The corresponding velocity fields read
$$
  \vv^{\gg_i} = \ff_i\ , \quad \vv^{\hh_{ij}} = \partial_i \vv^{\ff_j}
  + \partial_j \vv^{\ff_i}\ .
$$
Clearly, $\vv^{\gg_i}(\xi)$ has a Gaussian decay as $|\xi| \to \infty$, 
whereas $\vv^{\hh_{ij}}(\xi) \sim |\xi|^{-4}$. In particular, 
$\vv^{\hh_{ij}} \in \L^1(\real^3)$, but $\rho \vv^{\hh_{ij}} \notin 
\L^1(\real^3)$. 

\medskip
Now, we assume that $\ww \in \L^2(m)$ for some $m > 7/2$. Then $\ww$
can be decomposed as
\begin{equation}\label{wwdec}
  \ww(\xi) = \sum_{i=1}^3 \beta_i \ff_i(\xi) + 
  \sum_{i=1}^3 \gamma_i \gg_i(\xi) + \sum_{(ij)\in S} 
  \zeta_{ij}\hh_{ij}(\xi) + \tilde\ww(\xi)\ ,
\end{equation}
where the coefficients $\beta_i$ are defined in \reff{betadef} and 
$\gamma_i$, $\zeta_{ij}$ in \reff{gamzetdef}. The velocity field $\vv$ 
associated to $\ww$ has a similar decomposition:
\begin{equation}\label{vvdec}
  \vv(\xi) = \sum_{i=1}^3 \beta_i \vv^{\ff_i}(\xi) + 
  \sum_{i=1}^3 \gamma_i \vv^{\gg_i}(\xi) + \sum_{(ij)\in S} 
  \zeta_{ij}\vv^{\hh_{ij}}(\xi) + \tilde\vv(\xi)\ ,
\end{equation}
where $\tilde\vv$ is obtained from $\tilde\ww$ via the Biot-Savart law
\reff{BS3}. In view of \reff{betadef} and \reff{Bid}, it is clear that
$\beta_i = 0$ for $i = 1,2,3$ if and only if
\begin{equation}\label{Bcond}
   \int_{\real^3} \xi_i w_j(\xi)\d\xi = 0 \quad \hbox{for all }
   i,j \in \{1,2,3\}\ .
\end{equation}
On the other hand, it follows from \reff{gamzetdef} and \reff{Cid} that
$\zeta_{ij} = 0$ for all $(ij) \in S$ if and only if
\begin{equation}\label{Ccond}
\begin{array}{l}
  M^1_{22} = M^1_{33} = -2 M^2_{12} = -2 M^3_{13}\ , \cr
  M^2_{11} = M^2_{33} = -2 M^1_{12} = -2 M^3_{23}\ , \cr
  M^3_{11} = M^3_{22} = -2 M^1_{13} = -2 M^2_{23}\ ,
\end{array}
\qquad
\begin{array}{l}
  M^1_{11} = M^2_{22} = M^3_{33} = 0\ , \cr
  M^1_{23} = M^2_{13} = M^3_{12} = 0\ , \cr
\end{array}
\end{equation}
where
$$
   M^i_{jk} = \int_{\real^3} \xi_j \xi_k w_i(\xi)\d\xi = M^{i}_{kj}\ ,
   \quad i,j,k \in \{1,2,3\}\ .
$$
If in addition $\gamma_i = 0$ for $i = 1,2,3$, then $M^i_{jk} = 0$ for
all $i,j,k \in \{1,2,3\}$. 

\medskip
The main result of this section is the following estimate for the
velocity field in terms of the vorticity:

\begin{proposition}\label{velvort3} Let $\ww \in \L^2(m)$ for some 
$m \ge 0$, and let $\vv$ be the velocity field obtained from $\ww$ via 
the Biot-Savart law \reff{BS3}. Define the coefficients $\beta_i$ by 
\reff{betadef} and $\zeta_{ij}$ by \reff{gamzetdef}. Assume that either\\
\null~1) $0 \le m < 3/2$, or\\
\null~2) $3/2 < m < 5/2$, or\\
\null~3) $5/2 < m < 7/2$ and $\beta_i = 0$ for $i = 1,2,3$, or\\
\null~4) $7/2 < m < 9/2$, $\beta_i = 0$ for $i = 1,2,3$, and $\zeta_{ij} =0$
for $(ij) \in S$.\\
Then there exists $C > 0$ such that
\begin{equation}\label{3vel}
  |\rho^m \vv|_6 \le C|\rho^m \ww|_2\ ,
\end{equation}
where $\rho(\xi) = 1+|\xi|$. 
\end{proposition}

\bigskip\noindent{\bf Remark.} More generally, Proposition~\ref{velvort3} 
holds for any (not necessarily divergence free) vector field $\ww$ satisfying
$\rho^m \ww \in L^2(\real^3)^3$, provided that either\\
\null~1) $0 \le m < 3/2$, or\\
\null~2) $3/2 < m < 5/2$ and \reff{Aid} holds, or\\
\null~3) $5/2 < m < 7/2$ and \reff{Aid}, \reff{Bcond} hold, or\\
\null~4) $7/2 < m < 9/2$ and \reff{Aid}, \reff{Bcond}, \reff{Ccond} hold.\\
Roughly speaking, the result means that $\vv(\xi)$ decays like 
$|\xi|^{-m-\frac{1}{2}}$ as $|\xi| \to \infty$.

\medskip
The proof of Proposition~\ref{velvort3} is naturally divided into four 
steps. In the case 1), the bound \reff{3vel} is a direct consequence of 
\reff{SBS3} and of the following weighted Hardy-Littlewood-Sobolev 
inequality:

\begin{lemma}\label{L1}
If $-1/2 < m < 3/2$ and 
$$
  u(\xi) = \int_{\real^3} \frac{\omega(\eta)}{|\xi-\eta|^2}\d\eta\ ,
  \quad \xi \in \real^3\ ,
$$
then $|\rho^m u|_6 \le C|\rho^m\omega|_2$. 
\end{lemma}

\proof We use the dyadic decomposition
$$
   \real^3 = \bigcup_{j=0}^\infty B_j\ ,
$$
where $B_0 = \{\xi \in \real^3\,|\,|\xi| \le 1\}$ and $B_j = \{\xi \in 
\real^3\,|\, 2^{j-1} < |\xi| \le 2^j\}$ for $j \in \intplus^*$. Let 
$u_i = u\oone_{B_i}$ and $\omega_i = \omega\oone_{B_i}$, $i \in \intplus$. 
Clearly $u_i = \sum_{j \in \intplus} \Delta_{ij}$, where
$$
   \Delta_{ij}(\xi) = \oone_{B_i}(\xi) \int_{B_j}\frac{\omega_j(\eta)}
   {|\xi-\eta|^2}\d\eta\ .
$$
If $|i-j| \le 1$, it follows from \reff{HLSbis} that $|\Delta_{ij}|_6 
\le C|\omega_j|_2$. If $|i-j| \ge 2$, Young's inequality implies that 
$|\Delta_{ij}|_6 \le M_1^{3/4} M_2^{1/4}|\omega_j|_2$, 
where
$$
   M_1 = \sup_{\xi\in B_i}\Bigl(\int_{B_j} \frac{1}{|\xi-\eta|^3}\d\eta
   \Bigr)^{2/3}\ , \quad 
   M_2 = \sup_{\eta\in B_j}\Bigl(\int_{B_i} \frac{1}{|\xi-\eta|^3}\d\xi
   \Bigr)^{2/3}\ .
$$
If $i \ge j+2$, then $|\xi-\eta| \ge |\xi|-|\eta| \ge 2^{i-1}-2^j \ge 
2^{i-2}$ for all $\xi \in B_i$, $\eta \in B_j$. Thus $M_1 \le 
C2^{-2i}\mu(B_j)^{2/3} \le C2^{-2(i-j)}$ and $M_2 \le C2^{-2i}\mu(B_i)^{2/3}
\le C$ for some $C > 0$ independent of $i,j$, hence $|\Delta_{ij}|_6 \le 
C 2^{-\frac{3}{2}(i-j)}|\omega_j|_2$. If $j \ge i+2$, then $|\xi-\eta| 
\ge 2^{j-2}$ for all $\xi \in B_i$, $\eta \in B_j$, and a similar calculation 
shows that $|\Delta_{ij}|_6 \le C 2^{-\frac{1}{2}(j-i)}|\omega_j|_2$. 
Summarizing, we have shown that
$$
  |u_i|_6 \le C \sum_{j \in \intplus} K_{ij} |\omega_j|_2\ , \quad
  i \in \intplus\ ,
$$
where $K_{ij} = 2^{-|i-j|-\frac{1}{2}(i-j)}$. Now, by definition of the sets 
$B_i$, we have $|\rho^m u_i| \le C2^{mi}|u_i|$ and $|\rho^m \omega_j| \ge 
C 2^{mj}|\omega_j|$ for all $i,j \in \intplus$. It follows that 
$$
  |\rho^m u_i|_6 \le C\sum_{j \in \intplus} K_{ij}^{(m)} |\rho^m 
  \omega_j|_2\ , \quad i \in \intplus\ ,
$$
where $K_{ij}^{(m)} = 2^{-|i-j|+(m-\frac{1}{2})(i-j)}$. In particular, 
$|K_{ij}^{(m)}| \le 2^{-\alpha|i-j|}$ for some $\alpha > 0$, hence 
$K^{(m)}$ defines a bounded linear operator from $\ell^2(\intplus)$ into 
$\ell^6(\intplus)$. This concludes the proof. \QED

\medskip To prove \reff{3vel} in the remaining cases, we also need
the following variant of Lemma~\ref{L1}:

\begin{lemma}\label{L2}
If $-\frac{1}{2} < m < \frac{1}{2}$ and 
$$
  u(\xi) = \int_{\real^3} \frac{\omega(\eta)}{|\xi-\eta|}\d\eta\ ,
  \quad \xi \in \real^3\ ,
$$
then $|\rho^m u|_6 \le C|\rho^{m+1}\omega|_2$. 
\end{lemma}

\proof We use the same notations as in the preceding proof. 
If $|i-j| \le 1$, it follows from the Hardy-Littlewood-Sobolev 
inequality that $|\Delta_{ij}| \le C|\omega_j|_{6/5}$. By H\"older, 
$|\omega_j|_{6/5} \le C\mu(B_j)^{1/3}|\omega_j|_2 \le C 2^j |\omega_j|_2$,
hence $|\Delta_{ij}|_6 \le C 2^j |\omega_j|_2$. If $|i-j| \ge 2$, then 
$|\Delta_{ij}|_6 \le N_1^{3/4} N_2^{1/4} |\omega_j|_2$, where
$$
   N_1 = \sup_{\xi\in B_i}\Bigl(\int_{B_j} \frac{1}{|\xi-\eta|^{3/2}}\d\eta
   \Bigr)^{2/3}\ , \quad 
   N_2 = \sup_{\eta\in B_j}\Bigl(\int_{B_i} \frac{1}{|\xi-\eta|^{3/2}}\d\xi
   \Bigr)^{2/3}\ .
$$
Proceeding as above, we deduce that $|\Delta_{ij}|_6 \le C 2^{-\frac{1}{2}
|i-j|} 2^j |\omega_j|_2$ for all $i,j \in \intplus$. It follows that
$$
   |\rho^m u_i|_6 \le C \sum_{j \in \intplus}\tilde K_{ij}^{(m)} 
   |\rho^{m+1} \omega_j|_2\ ,
$$
where $\tilde K_{ij}^{(m)} = 2^{-\frac{1}{2}|i-j|+m(i-j)}$. Thus, 
if $|m| < 1/2$, $\tilde K^{(m)}$ defines a bounded linear operator from 
$\ell^2(\intplus)$ into $\ell^6(\intplus)$. \QED

\medskip We are now ready to prove Proposition~\ref{velvort3} in the
case 2). If $3/2 < m < 5/2$ and $\ww \in \L^2(m)$, then $\ww$ is 
integrable and \reff{Aid} holds. As a consequence, we can rewrite 
\reff{SBS3} in the form
\begin{equation}\label{I1}
  v_i(\xi) = -\frac{1}{4\pi} \sum_{j,k=1}^3 \epsilon_{ijk} 
  \int_{\real^3} \Bigl(\frac{\xi_j-\eta_j}{|\xi - \eta|^3} 
  - \frac{\xi_j}{|\xi|^3}\Bigr)w_k(\eta)\d\eta\ ,
\end{equation}
where $\epsilon_{ijk} = \sign(\sigma)$ if $(ijk)$ is a permutation 
$\sigma$ of $(123)$, and $\epsilon_{ijk} = 0$ otherwise. Using the identity
$$
  |\xi|^3 (\xi_j{-}\eta_j) -|\xi{-}\eta|^3 \xi_j = (\xi_j{-}\eta_j)|\xi|^2
  (|\xi|-|\xi{-}\eta|) + |\xi{-}\eta|(2\xi_j (\xi\cdot\eta) -\eta_j |\xi|^2 
  -\xi_j|\eta|^2)\ ,
$$
we obtain
\begin{eqnarray}\label{I2}
  \big| |\xi|^3 (\xi_j{-}\eta_j) -|\xi{-}\eta|^3 \xi_j \big| &\le& 
  C |\xi{-}\eta| |\xi| |\eta| (|\xi|+|\eta|) \\ \nonumber
  &\le& C (|\xi{-}\eta| |\xi|^2 |\eta| + |\xi{-}\eta|^2 |\xi| |\eta|)\ .
\end{eqnarray}
Thus, it follows from \reff{I1} that $|\vv(\xi)| \le C(f(\xi) +
g(\xi))$, where
$$
  f(\xi) = \frac{1}{|\xi|} \int_{\real^3} \frac{|\eta| |\ww(\eta)|}
    {|\xi-\eta|^2}\d\eta\ , \quad
  g(\xi) = \frac{1}{|\xi|^2} \int_{\real^3} \frac{|\eta| |\ww(\eta)|}
    {|\xi-\eta|}\d\eta\ .
$$
Since we already know that $|\vv|_6 \le C|\ww|_2$, it is sufficient
to bound $|\rho_0^m \vv|_6$, where $\rho_0(\xi) = |\xi|\oone_{\{|\xi|\ge 1\}}$.
Applying Lemma~\ref{L1} with $u(\xi) = |\xi|f(\xi)$ and $\omega(\eta) = 
|\eta||\ww(\eta)|$, we find $|\rho_0^m f|_6 \le |\rho^{m-1}u|_6 \le 
C|\rho^{m-1}\omega|_2 \le C|\rho^m \ww|_2$. Similarly, applying 
Lemma~\ref{L2} with $u(\xi) = |\xi|^2 g(\xi)$ and $\omega(\eta) = 
|\eta||\ww(\eta)|$, we find $|\rho_0^m g|_6 \le C|\rho^{m-2}u|_6 \le 
C|\rho^{m-1}\omega|_2 \le C|\rho^m \ww|_2$. Summarizing, we have shown 
that $|\rho^m \vv|_6 \le C|\rho^m \ww|_2$ if $3/2 < m < 5/2$. 

\medskip
We next assume that $5/2 < m < 7/2$ and that $\beta_i = 0$ for $i = 1,2,3$.
As was already observed, this implies that \reff{Bcond} holds. Thus 
\reff{I1} can be written in the form
\begin{equation}\label{I3}
  v_i(\xi) = -\frac{1}{4\pi} \sum_{j,k=1}^3 \epsilon_{ijk} 
  \int_{\real^3} \Bigl(\frac{\xi_j-\eta_j}{|\xi - \eta|^3} 
  - \frac{\xi_j-\eta_j}{|\xi|^3} -\frac{3(\xi\cdot\eta)\xi_j}{|\xi|^5}
  \Bigr)w_k(\eta)\d\eta\ .
\end{equation}
To bound the right-hand side, we observe that
\begin{equation}\label{I4}
\begin{array}{l}
  |\xi|^5 (\xi_j{-}\eta_j) -|\xi|^2 |\xi{-}\eta|^3(\xi_j{-}\eta_j)
   -3|\xi{-}\eta|^3 (\xi\cdot\eta)\xi_j \cr
  \hspace{2cm} = (\xi_j{-}\eta_j)|\xi|^2(|\xi|^3 - |\xi{-}\eta|^3 -3|\xi|
   (\xi\cdot\eta)) \cr
  \hspace{2cm} + \,3(\xi\cdot\eta)(|\xi|^3(\xi_j{-}\eta_j) - 
   |\xi-\eta|^3 \xi_j)\ .
\end{array}
\end{equation}
We claim that
\begin{equation}\label{I5}
  \big| |\xi|^3 - |\xi{-}\eta|^3 -3|\xi|(\xi\cdot\eta) \big| \le C
  |\eta|^2 (|\xi|+|\eta|)\ .
\end{equation}
Indeed, this follows from the identity
$$
\begin{array}{rcl}
  |\xi|^3 - |\xi{-}\eta|^3 -3|\xi|(\xi\cdot\eta) &=& 2(|\xi{-}\eta| -|\xi|)
    (\xi\cdot\eta) -|\xi{-}\eta||\eta|^2 \cr 
  &-& |\xi|(|\xi||\xi{-}\eta| + (\xi\cdot\eta)-|\xi|^2)\ ,
\end{array}
$$
and from the bound 
\begin{equation}\label{I5bis}
  0 \le |\xi||\xi{-}\eta| +(\xi\cdot\eta)-|\xi|^2 \le
  \frac{1}{2}|\eta|^2\ , 
\end{equation}
which is easily proved by setting $\eta = \xi-\zeta$, $\zeta \in \real^3$. 

Using \reff{I4} together with \reff{I2}, \reff{I5}, we obtain
$$
\begin{array}{l}
  \big| |\xi|^5 (\xi_j{-}\eta_j) -|\xi|^2 |\xi{-}\eta|^3(\xi_j{-}\eta_j)
   -3|\xi{-}\eta|^3 (\xi\cdot\eta)\xi_j \big| \cr
  \hspace{4cm} \le C |\xi{-}\eta| |\xi|^2 |\eta|^2 (|\xi|+|\eta|) \cr
  \hspace{4cm} \le C (|\xi{-}\eta||\xi|^3|\eta|^2 + |\xi{-}\eta|^2 |\xi|^2
  |\eta|^2)\ .
\end{array}
$$
Thus, it follows from \reff{I3} that $|\vv(\xi)| \le C(f(\xi) + g(\xi))$, 
where
$$
  f(\xi) = \frac{1}{|\xi|^2} \int_{\real^3} \frac{|\eta|^2 |\ww(\eta)|}
    {|\xi-\eta|^2}\d\eta\ , \quad
  g(\xi) = \frac{1}{|\xi|^3} \int_{\real^3} \frac{|\eta|^2 |\ww(\eta)|}
    {|\xi-\eta|}\d\eta\ .
$$
Applying Lemmas~\ref{L1} and \ref{L2} as in the previous case, we easily
obtain $|\rho_0^m f|_6 + |\rho_0^m g|_6 \le C|\rho^m \ww|_2$. This prove
\reff{3vel} in the case 3). 

\medskip Finally, we assume that $7/2 < m < 9/2$ and that $\beta_i = 0$
for $i = 1,2,3$, $\zeta_{ij} = 0$ for $(ij) \in S$. As we already 
remarked, this implies that \reff{Ccond} holds. Using \reff{Ccond},  
it is straightforward (but somewhat tedious) to verify that
\begin{equation}\label{J1}
  \sum_{j,k=1}^3 \epsilon_{ijk} \int_{\real^3} \Bigl(
  2|\xi|^2 (\xi\cdot\eta)\eta_j + \xi_j(|\xi|^2|\eta|^2 -5(\xi\cdot\eta)^2)
  \Bigr)w_k(\eta)\d\eta = 0\ , \quad i \in \intplus\ .
\end{equation}
Assuming \reff{J1}, we can rewrite \reff{I3} in the form
\begin{equation}\label{J2}
  v_i(\xi) = -\frac{1}{4\pi} \sum_{j,k=1}^3 \epsilon_{ijk} 
  \int_{\real^3} \frac{A_j(\xi,\eta)w_k(\eta)}{|\xi-\eta|^3 |\xi|^7}
  \d\eta\ ,
\end{equation}
where
\begin{eqnarray*}
  A_j(\xi,\eta) &=& |\xi|^7 (\xi_j{-}\eta_j) -|\xi|^4 |\xi{-}\eta|^3
  (\xi_j{-}\eta_j) \cr
  &-& 3|\xi|^2 |\xi{-}\eta|^3 (\xi\cdot\eta)(\xi_j{-}\eta_j) 
  -\frac{3}{2}\xi_j |\xi{-}\eta|^3 (5(\xi\cdot\eta)^2-|\xi|^2|\eta|^2)\ .
\end{eqnarray*}
We claim that
\begin{equation}\label{J3}
  |A_j| \le C|\xi{-}\eta| |\xi|^3 |\eta|^3 (|\xi|+|\eta|) \le 
  C (|\xi{-}\eta| |\xi|^4 |\eta|^3 + |\xi{-}\eta|^2 |\xi|^3 |\eta|^3)\ .
\end{equation}
Assuming \reff{J3} for the moment, we deduce from \reff{J2} that 
$|\vv(\xi)| \le C(f(\xi) + g(\xi))$, where
$$
  f(\xi) = \frac{1}{|\xi|^3} \int_{\real^3} \frac{|\eta|^3 |\ww(\eta)|}
    {|\xi-\eta|^2}\d\eta\ , \quad
  g(\xi) = \frac{1}{|\xi|^4} \int_{\real^3} \frac{|\eta|^3 |\ww(\eta)|}
    {|\xi-\eta|}\d\eta\ .
$$
Applying Lemmas~\ref{L1} and \ref{L2} again, we obtain $|\rho_0^m f|_6 + 
|\rho_0^m g|_6 \le C|\rho^m \ww|_2$. This proves \reff{3vel} in the 
case 4). 

It remains to establish \reff{J3}. We first remark that $A_j = (\xi_j{-}
\eta_j)|\xi|^2 B + C_j$, where
\begin{eqnarray*}
  B &=& |\xi|^5 -|\xi|^2 |\xi{-}\eta|^3 -3|\xi{-}\eta|^3(\xi\cdot\eta)
    -\frac{3}{2}|\xi| (5(\xi\cdot\eta)^2-|\xi|^2|\eta|^2)\ ,\cr
  C_j &=& \frac{3}{2}\bigl(|\xi|^3 (\xi_j{-}\eta_j) -\xi_j|\xi{-}\eta|^3
    \bigr) \bigl(5(\xi\cdot\eta)^2-|\xi|^2|\eta|^2\bigr)\ .
\end{eqnarray*}
In view of \reff{I2}, it is clear that $|C_j| \le C|\xi{-}\eta| |\xi|^3 
|\eta|^3 (|\xi|+|\eta|)$. Thus, it is sufficient to show that
\begin{equation}\label{J4}
  |B| \le C |\xi| |\eta|^3 (|\xi|+|\eta|)\ .
\end{equation}
To this end, we remark that
\begin{eqnarray}\label{J5}
  B &=& |\xi| D + (\xi\cdot\eta) (|\xi|^3 -|\xi{-}\eta|^3 -3|\xi|
    (\xi\cdot\eta)) \cr
  &+& (|\xi{-}\eta|-|\xi|) (4(\xi\cdot\eta)^2-|\xi|^2|\eta|^2) 
    -2|\eta|^2 (\xi\cdot\eta)|\xi{-}\eta|\ ,
\end{eqnarray}
where
\begin{eqnarray*}
  D &=& |\xi|^4 -|\xi|^2(\xi\cdot\eta) +\frac{1}{2}|\xi|^2|\eta|^2 
    -\frac{1}{2}(\xi\cdot\eta)^2 -|\xi|^3|\xi{-}\eta| \cr
  &=& \frac{1}{2} \bigl(|\xi||\xi{-}\eta| + (\xi\cdot\eta)-|\xi|^2\bigr)
    \bigl(|\xi|(|\xi{-}\eta|-|\xi|)-(\xi\cdot\eta))\ .
\end{eqnarray*}
Using \reff{I5bis}, we find $|D| \le \frac{1}{2}|\xi||\eta|^3$. 
Inserting this bound into \reff{J5} and using \reff{I5}, we obtain
\reff{J4}. This concludes the proof of Proposition~\ref{velvort3}. \QED

\medskip We conclude with two corollaries which are used in 
the preceding sections. 

\begin{corollary}\label{3velint}
Assume that $\ww \in \L^2(m)$ for some $m > \frac{5}{2}$, and denote by 
$\vv$ the velocity field obtained from $\ww$ via the Biot-Savart law 
\reff{BS3}. Then $\vv \in L^1(\real^3)^3$ if and only if $\beta_i = 0$ 
for $i = 1,2,3$, and in this case $\int_{\real^3}v_i(\xi)\d\xi = 0$ for 
$i = 1,2,3$. 
\end{corollary}

\proof Without loss of generality, we assume that $5/2 < m < 7/2$. 
For any $\ww \in \L^2(m)$, we have the decomposition
$$
   \ww = \sum_{i=1}^3 \beta_i \ff_i + \tilde \ww\ , \quad
   \vv = \sum_{i=1}^3 \beta_i \vv^{\ff_i} + \tilde \vv\ ,
$$
where the coefficients $\beta_i$ are defined in \reff{betadef}. Then
the remainder $\tilde \ww \in \L^2(m)$ fulfills the moment conditions 
\reff{Bcond}. By Proposition~\ref{velvort3}, the corresponding velocity
field satisfies $\rho^m \tilde \vv \in L^6(\real^3)^3$, hence $\tilde 
\vv \in \L^1(\real^3)$. On the other hand, it is easy to verify that
$\sum \beta_i \vv^{\ff_i} \in \L^1(\real^3)$ if and only if $\beta_i = 0$ 
for $i = 1,2,3$. Thus $\vv \in \L^1(\real^3)$ if and only if 
$\beta_1 = \beta_2 = \beta_3 = 0$. In this case, $\int_{\real^3} 
v_i(\xi)\d\xi = 0$ for $i = 1,2,3$ because $\div \vv =0$, see
\reff{Aid}. \QED

\begin{corollary}\label{3velintbis}
Assume that $\ww \in \L^2(m)$ for some $m > \frac{7}{2}$, and denote by 
$\vv$ the velocity field obtained from $\ww$ via the Biot-Savart law 
\reff{BS3}. Then $\rho \vv \in L^1(\real^3)^3$ if and only if 
$\beta_i = 0$ for $i = 1,2,3$ and $\zeta_{ij} = 0$ for $(ij) \in S$. 
In this case, the matrix $(b_{k\ell})$ defined by
\begin{equation}\label{bkldef}
   b_{k\ell} = \int_{\real^3} \xi_k v_\ell(\xi)\d\xi\ , \quad
   k,\ell \in \{1,2,3\}\ ,
\end{equation}
is skew-symmetric, and $b_{12} = \gamma_3$, $b_{23} = \gamma_1$, $b_{31} = 
\gamma_2$.
\end{corollary}

\proof Without loss of generality, we assume that $7/2 < m < 9/2$. 
If $\ww \in \L^2(m)$, then $\ww$ and $\vv$ can be decomposed according
to \reff{wwdec}, \reff{vvdec}. By construction, the remainder 
$\tilde \ww \in \L^2(m)$ has vanishing first-order and second-order moments. 
Applying Proposition~\ref{velvort3}, we deduce that the corresponding 
velocity field satisfies $\rho^m \tilde \vv \in L^6(\real^3)^3$, hence 
$\rho \tilde \vv \in L^1(\real^3)^3$. Therefore, using the expressions 
above of $\ff_i$, $\gg_i$, and $\hh_{ij}$, it is not difficult to
show that $\rho \vv \in L^1(\real^3)^3$ if and only if $\beta_i = 0$ 
for $i = 1,2,3$ and $\zeta_{ij} = 0$ for $(ij) \in S$. In this case, the 
matrix $(b_{k\ell})$ defined by \reff{bkldef} is skew-symmetric because
$\div \vv = 0$, see \reff{Bid}. 

Assume now that $\rho \vv \in L^1(\real^3)^3$, and consider the 
vector field $\AA(\xi) = \qq_3(\xi) \wedge \vv(\xi)$, where 
$\qq_3$ is defined in \reff{qqdef}. Then $\AA \in L^{3/2}(\real^3)^3$ 
and $\div \AA = \pp_3 \cdot \vv -\qq_3 \cdot \ww$, where $\pp_3 = 
\rot \qq_3$ is defined in \reff{ppdef}. It follows that
$$
   0 = \int_{\real^3} \div \AA \d\xi = \int_{\real^3} 
   (\pp_3 \cdot \vv -\qq_3 \cdot \ww)\d\xi = b_{12} - \gamma_3\ ,
$$
in view of \reff{ppdef}, \reff{bkldef}, and \reff{gamzetdef}. The relations
$b_{23} = \gamma_1$ and $b_{31} = \gamma_2$ are proved in a similar way. 
\QED

\bibliographystyle{plain}
\bibliography{ref}

\end{document}